\documentclass[10pt]{article}
\usepackage{amsmath, amsfonts, amsthm, amssymb, verbatim}
\usepackage{graphicx}
\usepackage[mathcal]{euscript}
\usepackage{amstext}
\usepackage{amssymb,mathrsfs}
\usepackage[ansinew]{inputenc}

\def\XXint#1#2#3{{\setbox0=\hbox{$#1{#2#3}{\int}$}
\vcenter{\hbox{$#2#3$}}\kern-.5\wd0}}

\newtheorem{theorem}{Theorem}[section]
\newtheorem{proposition}[theorem]{Proposition}
\newtheorem{remark}[theorem]{Remark}
\newtheorem{lemma}[theorem]{Lemma}

\newtheorem{definition}[theorem]{Definition}

\input{tcilatex}
\begin{document}

\title{The generalized Buckley-Leverett System. Solvability}
\author{Nikolai Chemetov$^1$, Wladimir Neves$^2$}
\date{}
\maketitle

\begin{abstract}
We propose a new approach to the mathematical modeling of the Buckley-
Leverett system, which describes two-phase flows in porous media.
Considering the initial-boundary value problem for a deduced model, we prove
the solvability of the problem. The solvability result relies mostly on the
kinetic method.
\end{abstract}

\tableofcontents

\footnotetext[1]{%
Centro de Matemática e Aplicações Fundamentais, Universidade de Lisboa, Av.
Prof. Gama Pinto 2, 1649-003 Lisboa, Portugal. E-mail: \textsl{%
chemetov@ptmat.fc.ul.pt.}} \footnotetext[2]{%
Institute of Mathematics, Federal University of Rio de Janeiro, C.P. 68530,
Cidade Universitária, 21945-970, Rio de Janeiro, Brazil. E-mail: \textsl{%
wladimir@im.ufrj.br.} \newline
\textit{To appear in:}
}

%
%

%



\section{Introduction}

\label{s.1}

We propose a new and more general approach to the mathematical study of the
Buckley-Leverett system, which describes two-phase flows in porous media,
that is, the simultaneous motion~of two immiscible incompressible liquids
(e.g. water and oil) in a porous medium. This study is of practical interest
in connection with the planning and operation of oil wells and also brings
some challenging mathematical questions. Indeed, the mathematical standard
model of the Buckley-Leverett system is given by a scalar multidimensional
conservation law, which describes the evolution of the saturation according
to the seepage velocity field, where this field is given by the Darcy's Law
(empirical) equation, i.e.
\begin{equation*}
\text{seepage velocity= equivalent mobility(saturation) $\times $ pressure
gradient}.
\end{equation*}%
Since the equivalent mobility of the porous medium is a function of the
saturation, which should be a bounded and measurable function, this brings
an enormous mathematical difficult. We have to solve a scalar conservation
laws in the class of roughly coefficients. Even nowadays with the best
results of Panov for scalar conservation laws with discontinuous flux
published in ARMA, see \cite{EP}, we are not allowed to solve the
Buckley-Leverett system in this way. In order to pass the above difficult,
the Buckley-Leverett system has been significantly simplified in many works,
for instance see: \ \ Córdoba, Gancedo, Orive \cite{DCFGRO}, Frid \cite{Frid}%
, \ Perepelitsa, Shelukhin \cite{MPVS}. In the article of Luckhaus,
Plotinikov \cite{SLPIP} has considered a stationary case of the
Buckley-Leverett system. Many authors have proposed interesting ideas, but
most of them focused on the saturation equation. Mainly reducing the
Buckley-Leverett system to a (non)degenerate elliptic-parabolic partial
differential system, here we address some of the important works on this
subject:\ Antontsev, Kazikhov, Monakhov \cite{ant}, Arbogast \ \cite{arb},
Chen \cite{chen}, Lenzinger, Schweizer \cite{len}, Sazenkov \cite{SAS} and
further references cited therein.

\smallskip In the present work we change the focus and put more attention to
the equation of velocity. So we propose a generalized Darcy's law equation,
which is no physically longer than the standard one. This new formulation
bring to us enough regularity of the seepage velocity field in order to
obtain solvability of the system using the nice idea of Kinetic Theory. This
is the most part of the motivation to introduce a general approach to the
mathematical study of the Buckley-Leverett system. In the rest of the
introduction, we give a general presentation of the generalized Darcy's law
equation in a homogeneous and isotropic medium for one phase flow.

The theory of flows in porous media has also a number of similarities with
several other processes occurring in the continuum physics as, for instance,
problems of infiltration, displacement of electricity through dielectric
media, heat transfer, etc. Indeed, after suitable averaging the porous media
and the liquids filling them can be regarded as continuous medium under
natural assumptions made about the pore system, see Scheidegger \cite%
{AES,AES1}. Analogously to the resistance for the conductors of electricity,
we have here the porosity as a characteristic of the porous media.

One observes that, for very short time scales or high frequency
oscillations, a time derivative of flux may be added to Darcy's law, which
results in valid solutions at very small times. In heat transfer, a similar
idea is called the Cattaneo's law which is a modified version of the
standard Fourier one, hence in analogy we have the following equation for
the velocity field $\mathbf{v}$,
\begin{equation}
\tau \;\partial _{t}\mathbf{v}+\frac{\mu }{\kappa }\mathbf{v}=-\nabla p,
\label{1}
\end{equation}%
where $p$ is the pressure, $\mu $ is the dynamic viscosity, $\kappa $ is the
permeability and $\tau $ is a very small time constant. The parameter $\tau $
causes this equation to reduce to the normal form of Darcy's law at usual
times. The main reason for doing this is that the regular groundwater flow
equation (diffusion equation) leads to singularities at constant head
boundaries at very small times. Analogously to the heat transfer case, this
form is more mathematically rigorous, which leads to a hyperbolic
groundwater flow equation.

Another extension to the traditional form of Darcy's law is Brinkman´s term
(introduced in 1947 by Brinkman), which is used to account for transitional
flow between boundaries,
\begin{equation}
\nu \,\Delta \mathbf{v}+\frac{\mu }{\kappa }\mathbf{v}=-\nabla p,
\label{222}
\end{equation}%
where $\nu $ is an effective viscosity term. This correction term accounts
for flow through medium where the grains of the media are porous themselves.
In this paper, we are going to consider both extensions \eqref{1} and %
\eqref{222} of Darcy´s law equation for positive viscosity $\nu $ and
non-negative parameter $\tau $. The combination of \eqref{1} and \eqref{222}
is sometimes called Brinkman-Forchheimer equation in a porous media
literature, see for instance \cite{sheu}, \cite{straughan} and \cite{wang}.

\smallskip It is important to observe that, such generalized Darcy's law
models are described also considering the homogenization theory, see for
instance \cite{hornung}, and we address further \cite{straughan}.


\bigskip

\section{Non-linear porous-media theory}

\label{SNLPMT}

In this section we are going to formulate the porous-media theory for two
immiscible incompressible liquids in a porous medium. Here, as it is
standard in the formalism of continuum mechanics, we have enough regularity,
integrability, etc. of the involved functions to proceed any mathematical
computation. Moreover, we assume some simplifications and analogies to
obtain our model, but taken out any un-physical considerations.

Let $(t,\mathbf{x})\in \mathbb{R}\times \mathbb{R}^{d}$ $(d=1,2,3)$ be the
points in the time-space domain. First, let us consider the porosity $m(%
\mathbf{x})$, which is the proportion of the pore volume in an infinitesimal
part of the porous medium containing the point $\mathbf{x}\in \Omega $.

One could describes the problem of two-phase flow in a porous medium with
the two main elements $\mathbf{v}_{i}(t,\mathbf{x})$, $s_{i}(t,\mathbf{x})$ $%
(i=1,2)$ respectively the velocity field of each fluid, which takes value in
$\mathbb{R}^{d}$, and the saturation of each component, which is a scalar
function. In fact, the saturation of each component represents the local
proportion of the pore space occupied by the $i^{th}$-phase, thus we must
have $0\leq s_{i}\leq 1$ $(i=1,2)$ and
\begin{equation}
s_{1}+s_{2}=1.  \label{SE}
\end{equation}%
Moreover, the velocity field of each component is obtained from an average
of the flow rate of the $i^{th}$-phase divide by an unitary area, used
called seepage velocity. Related to the speed of the velocities, we are not
going to consider the non-linear convection terms.


\bigskip Now, we are in position to present the equations concerned the
immiscible incompressible multiple phase flow. At this point, we follow
reference \cite{AES1} Part IX, and address also \cite{AES}. The evolution of
the saturation is driven by the velocity field described by the following
continuity equation
\begin{equation}
m\;\partial _{t}(\rho _{i}\,s_{i})+\mathrm{div}(\rho _{i}\,\mathbf{v}%
_{i})=0,\qquad (i=1,2)  \label{SEC}
\end{equation}%
where $\rho _{i}$ is mass density of the $i^{th}$-phase of the porous medium
and the velocity field satisfies the generalized Darcy's law equation
\begin{equation}
\begin{aligned} \tau_i \, \rho_i \, \partial_t \mathbf{v}_i &- \nu_i \,
\Delta \mathbf{v}_i + \frac{\mu_i}{k_0 \, k_{ri}(s_1)} \mathbf{v}_i &= -
\nabla p_i + \rho_i \, g \, h, \qquad(i=1,2) \end{aligned}  \label{VEC}
\end{equation}%
where for each component $i=1,2$, $p_{i}(t,\mathbf{x})$ is the pressure, $%
\nu _{i}$, $\mu _{i}$ respectively the viscosity and dynamic viscosity, $%
k_{ri}$ is the relative permeability and $\tau _{i}$ is the time-delay
parameter. Moreover, $k_{0}(\mathbf{x})$ is the absolute permeability of the
porous medium and $\rho \,g\,h$ is the external gravitational force, which
could be dropped, since we are considering an horizontal domain (reservoir)
whose height is negligible compared to the other dimensions. As it stands,
the equations \eqref{SEC} and \eqref{VEC} form a system of four partial
differential equations, where the unknowns are the velocity and saturation
of each component. The pressure is obtained a posteriori by the velocity
field, exactly as a Lagrangian multiplier in the Leray's theory for
incompressible Navier-Stokes equations.

\medskip As it is standard in porous media theory, in order to simplify the
model, we assume that $m=k_{0}=\rho _{1}=\rho _{2}\equiv 1$, and further
\begin{equation*}
\tau _{1}=\tau _{2}=:\tau ,\quad \nu _{1}=\nu _{2}=:\nu .
\end{equation*}%
In fact, the time-delay $\tau $ is a nonnegative very-small parameter, and
here, we are going to consider two cases, that is $\tau >0$ and $\tau =0$.
The viscosity positive parameter $\nu $ is also very small, i.e. $0<\nu
<\!<1 $. Moreover, by the Laplace's formula (experimental one), it follows
that
\begin{equation*}
p_{1}(t,\mathbf{x})-p_{2}(t,\mathbf{x})=p_{c}(\mathbf{x},s_{1}),
\end{equation*}%
where $p_{c}$ is the capillarity pressure, and by the Buckley-Leverett
assumptions, we suppose that $p_{c}\equiv 0$, thus we have $p_{1}=p_{2}=:p$.
Then, from \eqref{SEC} and \eqref{VEC} we obtain respectively
\begin{equation}
\hspace{25pt}\partial _{t}s_{i}+\mathrm{div}\left( \mathbf{v}_{i}\right) =0
\label{SECS}
\end{equation}%
and
\begin{equation}
\begin{aligned} \tau \, \partial_t \mathbf{v}_i &- \nu \, \Delta_x
\mathbf{v}_i + \lambda_i \mathbf{v}_i &= - \nabla_x p, \end{aligned}
\label{VECS}
\end{equation}%
where $\lambda _{i}=\mu _{i}/k_{ri}$ $(i=1,2)$.

\medskip Now, we are going to proceed in order to obtain the final model,
which is written in terms of the (total) velocity $\mathbf{v}=\mathbf{v}_{1}+%
\mathbf{v}_{2}$. We derive it assuming temporarily the one-dimensional case
and denote $F(t,x)=-\partial _{x}p(t,x)$ for simplicity. First, let us
understand precisely the rule of the parameter $\tau $ in the generalized
Darcy's law equation. Set $t=\tau \,\hat{t}$ and define
\begin{equation*}
\mathbf{v}(t,x)=:\hat{\mathbf{v}}(\hat{t},x),\quad F(t,x)=:\hat{F}(\hat{t}%
,x),
\end{equation*}%
where we have dropped the under-script $(i=1,2)$ in order to simplify the
notation. Therefore, we have for very short times, i.e. $t=O(\tau )$, that
the velocity field $\mathbf{v}(t,x)$ behaves like in "normal" times, since $%
\mathbf{v}(t,x)=\hat{\mathbf{v}}(O(1),x)$, where $\hat{\mathbf{v}}$
satisfies the equation
\begin{equation*}
\partial _{\hat{t}}\hat{\mathbf{v}}-\nu \,\partial _{x}^{2}\hat{\mathbf{v}}%
+\lambda \,\hat{\mathbf{v}}=\hat{F}.
\end{equation*}%
On the other hand, for $t=O(1)$, $\mathbf{v}(t,x)$ behaves like in permanent
regime, since $\mathbf{v}(t,x)=\hat{\mathbf{v}}(O(1/\tau ),x)$.
Consequently, for each $\delta >0$ we are allowed to suppose that $\mathbf{v}%
(t+\delta ,x)\simeq \mathbf{v}(t,x)$ for usual times $t=O(1)$. Further, we
apply formally the Faedo-Galerkin's method to equation \eqref{VECS}, i.e.,
we consider
\begin{equation*}
\mathbf{v}(t,x)=\sum_{n=-\infty }^{\infty }\mathbf{v}_{n}(t)\exp (inx)\quad
\text{and}\quad F(t,x)=\sum_{n=-\infty }^{\infty }c_{n}(t)\exp (inx).
\end{equation*}%
Then, from \eqref{VECS} we obtain for each $n$
\begin{equation*}
\tau \,\mathbf{v}_{n}^{\prime }(t)+\Lambda _{n}\mathbf{v}_{n}(t)=c_{n}(t),
\end{equation*}%
where $\Lambda _{n}:=\lambda +n^{2}\nu $. Let $\delta >0$ be sufficiently
small, we resolve the above differential equation from $t$ to $t+\delta $
for some usual time $0<t=O(1)$, that is
\begin{equation*}
\mathbf{v}_{n}(t)\;\Big(e^{(\Lambda _{n}/\tau )\;(t+\delta )}-e^{(\Lambda
_{n}/\tau )\;t}\Big)=\frac{1}{\tau }\int_{t}^{t+\delta }c_{n}(\xi
)\;e^{(\Lambda _{n}/\tau )\;\xi }\;d\xi ,
\end{equation*}%
where we have used $\mathbf{v}_{n}(t+\delta )=\mathbf{v}_{n}(t)$. Hence
dividing by $\delta $ and taking the limit as $\delta \rightarrow 0^{+}$, we
have
\begin{equation*}
\mathbf{v}_{n}(t)=\frac{1}{\Lambda _{n}}\;c_{n}(t).
\end{equation*}%
Therefore, it follows that
\begin{equation*}
\begin{aligned} \mathbf{v}(t,x)&= \sum_{n=-\infty}^\infty
\frac{1}{\Lambda_{n}} \; c_n(t) \exp(i n x) \\ &= \frac{1}{\Lambda} \;
F(t,x), \end{aligned}
\end{equation*}%
where $1/\Lambda $ is the value of the absolutely convergent series $\sum
1/(\lambda +n^{2}\nu )$. From the above expression, i.e. $\Lambda _{i}\,%
\mathbf{v}_{i}=F$ $(i=1,2)$, we obtain
\begin{equation*}
\mathbf{v}:=\mathbf{v}_{1}+\mathbf{v}_{2}=\Big(\frac{\Lambda _{1}+\Lambda
_{2}}{\Lambda _{1}\;\Lambda _{2}}\Big)F(t,x),
\end{equation*}%
or also denoting
\begin{equation*}
\Lambda _{eq}:=\frac{\Lambda _{1}\;\Lambda _{2}}{\Lambda _{1}+\Lambda _{2}},
\end{equation*}%
we have
\begin{equation*}
F=\Lambda _{eq}\;\mathbf{v}=\Lambda _{1}\;\mathbf{v}_{1}=\Lambda _{2}\;%
\mathbf{v}_{2},
\end{equation*}%
that is,
\begin{equation*}
\mathbf{v}_{i}(t,x)=\frac{\Lambda _{eq}}{\Lambda _{i}}\;\mathbf{v}%
(t,x)\qquad (i=1,2).
\end{equation*}%
Finally, taking as motivation the above formulation, we derive our
porous-media generalized model for two immiscible incompressible liquids in
a porous media in the following way. From equation \eqref{SECS} written for
each component and added, we obtain
\begin{equation}
\mathrm{div}\mathbf{v}=0.  \label{DVF}
\end{equation}%
Moreover, denoting $u:=s_{1}$, we have
\begin{equation}
\partial _{t}u+\mathbf{v}\cdot \nabla g(u)=0,  \label{EQ100}
\end{equation}%
where $g(u):=\Lambda _{eq}/\Lambda _{1}$. Finally, taking in account the
parabolic/eliptic equation \eqref{VECS}, we have
\begin{equation*}
\tau \;\partial _{t}\mathbf{v}_{1}-\nu \Delta \mathbf{v}_{1}+\lambda _{1}%
\mathbf{v}_{1}=-\nabla p,
\end{equation*}%
and
\begin{equation*}
\tau \;\partial _{t}\mathbf{v}_{2}-\nu \Delta \mathbf{v}_{2}+\lambda _{2}%
\mathbf{v}_{2}=-\nabla p,
\end{equation*}%
hence adding these two equations, we get
\begin{equation}
\tau \;\partial _{t}\mathbf{v}-\nu \Delta \mathbf{v}+h(u)\mathbf{v}=-2\nabla
p,  \label{VFF}
\end{equation}%
where
\begin{equation*}
\begin{aligned} h(u)&= \Big( \frac{\lambda_1}{\Lambda_1}+
\frac{\lambda_2}{\Lambda_2}\Big) \Lambda_{eq}= \frac{\lambda_1 \Lambda_2 +
\lambda_2 \Lambda_1} {(\Lambda_1 + \Lambda_2)}. \end{aligned}
\end{equation*}

The above deduced model, that is, equations \eqref{DVF}, \eqref{EQ100} and %
\eqref{VFF} will be called \textbf{Stokes-Buckley-Leverett system} (or for
brevity \textbf{Stokes }\textit{B-L} system) for $\tau \neq 0$. Moreover,
when $\tau= 0$, we are going to say \textbf{quasi-stationary} \textbf{Stokes
}\textit{B-L} system.



\vspace{1pt}

\section{Functional notation and background}

\label{SSFNB}

Let $T>0$ be any fixed real number and $\Omega \subset \mathbb{R}^{d}$ (with
$d=1,2$ or $3)$ is an open and bounded domain having a $C^{2}-$ smooth
boundary $\Gamma .$ We define by
\begin{equation*}
\Omega _{T}:=\Omega \times (0,T),\qquad \Gamma _{T}:=\Gamma \times (0,T).
\end{equation*}%
The outside normal to $\Omega $ at $\mathbf{x\in }\Gamma $ is denoted by $%
\mathbf{n=n}(\mathbf{x}).$ \

In the paper we will use the standard notations for the Lebesgue function
space $L^{p}(\Omega )$ and the Sobolev spaces $W^{s,p}(\Omega )$ and $%
H^{s}(\Omega )\equiv W^{s,2}(\Omega )$ where a real $p\geqslant 1$ is the
integrability indice and a real $s\geqslant 0$ is the smoothness indice. The
vector counterparts of these spaces are denoted by $\mathbf{L}^{2}(\Omega
)=(L^{2}(\Omega ))^{d},$ \ $\mathbf{W}^{s,p}(\Omega ):=\left( W^{s,p}(\Omega
)\right) ^{d}$\ and $\mathbf{H}^{s}(\Omega ):=(H^{s}(\Omega ))^{d}.$ Let us
point that by Theorem 1.2 of \ \cite{Temam} for any $\mathbf{u}\in \mathbf{L}%
^{2}(\Omega ),$ satisfying $\mathrm{div}\left( \mathbf{u}\right) =0$\quad in
$\mathcal{D}^{^{\prime }}(\Omega ),$\quad the normal component of $\mathbf{u,%
}$ i.e. $\mathbf{u}_{\mathbf{n}}\mathbf{:=u}\cdot \mathbf{n,}$ exists and
belongs to $H^{-1/2}(\Gamma ).$ We will also use the following divergence
free spaces

\begin{eqnarray*}
\mathbf{V}^{s}(\Omega ) &:&=\{\mathbf{u}\in \mathbf{H}^{s}(\Omega ):\quad
\mathrm{div}\left( \mathbf{u}\right) =0\quad \text{in}\quad \mathcal{D}%
^{^{\prime }}(\Omega ),\quad \int_{\Gamma }\mathbf{u}_{\mathbf{n}}\ d\mathbf{%
x}=0\}, \\
\mathbf{V}^{s}(\Gamma ) &:&=\{\mathbf{u}\in \mathbf{H}^{s}(\Gamma ):\quad
\int_{\Gamma }\mathbf{u}_{\mathbf{n}}\ d\mathbf{x}=0\},\quad \quad \mathbf{V}%
^{-s}(\Gamma ):=\left( \mathbf{V}^{s}(\Gamma )\right) ^{\prime }
\end{eqnarray*}%
and

\begin{equation*}
\mathbf{G}(\Gamma _{T}):=\left\{ \mathbf{u}\in L^{2}(0,T;\mathbf{V}%
^{1/2}(\Gamma )):\quad ~\partial _{t}\mathbf{u}\in L^{2}(0,T;\mathbf{V}%
^{-1/2}(\Gamma )\right\} .
\end{equation*}

Let us formulate some results for the Stokes equations, used in the sequel.
We consider the system
\begin{equation}
\begin{cases}
-\nu \Delta \mathbf{v}=-\nabla p,\qquad \mathrm{div}(\mathbf{v})=0\quad
\text{in }\Omega , \\
\mathbf{v}=\mathbf{b}\quad \text{on }\Gamma .%
\end{cases}
\label{eq2}
\end{equation}%
The proof of the following result has been done by Cattabriga in \cite{cat}
(see also Theorem 3 with Remarks 2 of \cite{galdi}).

\begin{proposition}
\label{prop2} If $\;\mathbf{b}\in \mathbf{H}^{s-1/2}(\Gamma )$ for $s=0$ or $%
s=1,$ then there exists an unique weak solution $\mathbf{v}\in \mathbf{V}%
^{s}(\Omega )$ of \eqref{eq2}, such that
\begin{equation*}
\Vert \mathbf{v}\Vert _{\mathbf{V}^{s}(\Omega )}\leqslant C\Vert \mathbf{b}%
\Vert _{\mathbf{H}^{s-1/2}(\Gamma )}.
\end{equation*}
\end{proposition}

\vspace{1pt}

\section{Statement of the \textbf{Stokes }\textit{B-L} system}

\label{SSOP}

\vspace{1pt}In this section we are going to formulate the mathematical
problem. Let us assume that $\tau ,\nu $ are given positive fixed
parameters. We are concerned with the following initial-boundary value
problem, denoted as \textbf{IBVP}$_{\tau }$:

\textit{Find a pair }$(u,\mathbf{v})=(u(t,\mathbf{x}),\mathbf{v}(t,\mathbf{x}%
)):\Omega _{T}\rightarrow \mathbb{R}\times \mathbb{R}^{d}$\textit{\ solution
to the Stokes-Buckley-Leverett system in the domain } $\Omega _{T}$
\begin{eqnarray}  \label{EDFV}
\partial _{t}u+\mathrm{div}\big(\mathbf{v}\ g(u)\big) &=&0,  \label{ETU} \\
\tau \,\partial _{t}\mathbf{v}-\nu \Delta \mathbf{v}+h(u)\mathbf{v}
&=&-\nabla p,\qquad \mathrm{div}(\mathbf{v})=0,  \label{EPEV}
\end{eqnarray}%
\textit{satisfying the boundary conditions}%
\begin{equation}
(u,\mathbf{v})=(u_{b},\mathbf{b})\quad \text{on }\Gamma _{T},  \label{BC}
\end{equation}%
\textit{\ and the initial conditions }%
\begin{equation}
(u,\mathbf{v})|_{t=0}=(u_{0},\mathbf{v}_{0})\quad \text{in }\Omega .
\label{IC}
\end{equation}

We assume that our data satisfy the following regularity properties%
\begin{eqnarray}
g,\ h &\in &W_{\,\mathrm{loc}}^{1,\infty }(\text{$\mathbb{R}$})\qquad \text{%
with}\quad 0 < h_0 \leq h(u),\text{ }  \notag \\
0 &\leqslant &u_{b}\leqslant 1\quad \text{on }\Gamma _{T},  \notag \\
0 &\leqslant &u_{0}\leqslant 1\quad \text{on }\Omega ,  \label{regular1}
\end{eqnarray}

and
\begin{eqnarray}
\mathbf{v}_{0} &\in &\mathbf{V}^{0}(\Omega )\quad \text{and}\quad \mathbf{b}%
\in \mathbf{G}(\Gamma _{T}),\quad \text{such that}  \notag \\
\mathbf{b}(0)\cdot \mathbf{n} &=&\mathbf{v}_{0}\cdot \mathbf{n}\quad \quad
\text{in}\quad H^{-1/2}(\Gamma ).  \label{regular2}
\end{eqnarray}

\bigskip

Now, since equation \eqref{ETU} is a hyperbolic scalar conservation law, the
saturation function $u$ may admit shocks. Therefore, in order to select the
more correct physical solution, we need the entropy concept as given at the
following

\begin{definition}
\label{DEF} A pair $\mathbf{F}(u):=(\eta (u),q(u))$ is called an entropy
pair for \eqref{ETU}, if $\eta :\mathbb{R}\rightarrow \mathbb{R}$ is a
Lipschitz continuous and also convex function and the function $q:\mathbb{R}%
\rightarrow \mathbb{R}$ satisfies%
\begin{equation}
q^{\prime }(u)=\eta ^{\prime }(u)\;g^{\prime }(u)\quad \text{for a. a. }u\in
\mathbb{R}.  \label{DEFE}
\end{equation}
\end{definition}

Certainly, the most important example of entropy pairs are given by the Kru%
\u{z}kov's entropies. Here, we consider the following parameterized family
of Kru\u{z}kov's entropy pairs for \eqref{ETU}%
\begin{equation}
\mathbf{F}(u,v)=\Big(|u-v|,\,\mathrm{sgn}(u-v)\big(g(u)-g(v)\big)\Big)
\label{KE}
\end{equation}%
for each $v\in \mathbb{R}.$ We remark that any smooth entropy pair $\mathbf{F%
}(u):=(\eta (u),q(u))$ for \eqref{ETU} can be recovered by the family given
by \eqref{KE}. The inverse one is also true, i.e. any entropy pair $\mathbf{F%
}(u):=(\eta (u),q(u))$ given by \eqref{KE} can be recovered by a family of
smooth entropy pairs. \ In fact, this result follows for any entropy by a
standard regularization argument.

Another two examples of parameterized family of entropy pairs for \eqref{ETU}
are
\begin{equation}
\mathbf{F}^{\pm }(u,v)=\Big(~|u-v|^{\pm },~\mathrm{sgn}^{\pm }\left(
u-v\right) \left( g(u)-g(v)\right) ~\Big)  \label{KE1}
\end{equation}%
for each $v\in \mathbb{R},$ which will be useful in the Kinetic formulation
(see Section \ref{SSKL}). Here
\begin{equation*}
\mathrm{sgn}^{+}\left( v\right) :=\left\{
\begin{array}{l}
1,\quad \text{if }v>0, \\
0,\quad \text{if }v\leqslant 0,%
\end{array}%
\right. \quad \quad \mathrm{sgn}^{-}\left( v\right) :=\left\{
\begin{array}{l}
0,\quad \text{if }v<0, \\
1,\quad \text{if }v\geqslant 0%
\end{array}%
\right.
\end{equation*}%
and $\ \ \ |v|^{\pm }:= \max\left\{ \pm v,~0\right\} ,$ respectively.

The following definition tells us in which sense a pair of functions $(u,%
\mathbf{v})$ is a weak solution of \textbf{IBVP}$_{\tau }$: \eqref{ETU}-%
\eqref{IC}.

\begin{definition}
\label{DGS11} A pair of functions
\begin{equation*}
(u,\mathbf{v})\in L^{\infty }(\Omega _{T})\times L^{2}(0,T;\mathbf{V}%
^{1}(\Omega ))
\end{equation*}%
is called a weak solution to the \textbf{IBVP}$_{\tau }$, if the pair $(u,%
\mathbf{v})$ satisfies the integral inequality
\begin{align}
\iint_{\Omega _{T}}& \big(|u-v|\;\phi _{t}+\,\mathrm{sgn}(u-v)\big(g(u)-g(v)%
\big)\;\mathbf{v}\cdot \nabla \phi \big)\ d\mathbf{x}\,dt  \notag \\
& +\int_{\Gamma _{T}}M\;|u_{b}-v|\phi \ d\mathbf{x}\,dt+\int_{\Omega
}|u_{0}-v|\phi (0,x)\ d\mathbf{x}\geq 0  \label{DGS}
\end{align}%
for any fixed $v\in \mathbb{R},$ where $M:=\emph{K}|\mathbf{b}_{\mathbf{n}}%
\mathbf{|}$ defined on $\Gamma _{T}$ with $\ \emph{K:=}||g^{\prime
}||_{L^{\infty }(\mathbb{R})}$\ and for any nonnegative function $\phi \in
C_{0}^{\infty }((-\infty ,T)\times \mathbb{R}^{d})$ and also the following
integral identity
\begin{equation}
\int_{\Omega _{T}}\left[ \tau \ \mathbf{v}\cdot \boldsymbol{\psi }_{t}-\nu
\,\nabla \mathbf{v}:\nabla \boldsymbol{\psi }-h(u)\ \mathbf{v}\cdot
\boldsymbol{\psi }\right] \,d\mathbf{x}dt+\tau \int_{\Omega }\mathbf{v}%
_{0}\cdot \boldsymbol{\psi }(0)\,d\mathbf{x}=0  \label{V}
\end{equation}%
holds for any $\boldsymbol{\psi }\in \mathbf{C}^{1}(\overline{\Omega }_{T}),$
such that $\boldsymbol{\psi }=0$ at $t=T$ and on $\Gamma _{T}.$ Moreover the
trace of $\mathbf{v}$ is equal to $\mathbf{b}$ on $\Omega .$ Here $\ \nabla
\mathbf{v}:\nabla \boldsymbol{\psi }=\sum_{j=1}^{d}\partial _{x_{j}}\mathbf{v%
}\cdot \partial _{x_{j}}\boldsymbol{\psi}.$
\end{definition}

For more complete discussions on this concept of weak entropy solutions for
hyperbolic conservation law \eqref{ETU} (with boundary conditions), we refer
to Otto \cite{O}, Neves \cite{WN1}, Chen, Frid \cite{CF} (see Theorem 4.1)
and Malek et all \cite{malek} (see Lemma 7.24 and Theorem 7.31),\textbf{\ }
further the Dafermos' treatise book \cite{Dafermos}.

\begin{theorem}
\label{ETGS} If the data $g,$ $h,$ $u_{b},$ $u_{0},$ $\mathbf{v}_{0},$ $%
\mathbf{b}$ fulfills the regularity properties \eqref{regular1}-%
\eqref{regular2}, then the \textbf{IBVP}$_{\tau }$ \ has a weak solution
\begin{equation*}
(u,\mathbf{v})\in L^{\infty }(\Omega _{T})\times L^{2}(0,T;\mathbf{V}%
^{1}(\Omega ))\cap H^{1}(0,T;\mathbf{V}^{-1}(\Omega )),
\end{equation*}%
satisfying%
\begin{eqnarray}
&&0\leqslant u\leqslant 1\text{\qquad a . e. in \quad }\Omega _{T},  \notag
\\
&\Vert \sqrt{\tau }\mathbf{v}\Vert _{C([0,T];\mathbf{V}^{0}(\Omega
))}&\!\!\!\!\!+\Vert \mathbf{v}\Vert _{L^{2}(0,T;\mathbf{V}^{1}(\Omega
))}+\tau \Vert \mathbf{v}\Vert _{H^{1}(0,T;\mathbf{V}^{-1}(\Omega ))}\leq C,
\label{r}
\end{eqnarray}%
where $C$ is a positive constant independent of $\tau .$
\end{theorem}

\section{Existence of weak solution}

\label{SEGS}

\bigskip

\subsection{Parabolic approximation}

\label{SSPA}

In order to show the existence of a weak solution for the \textbf{IBVP}$%
_{\tau }$, \ first we study the following approximated parabolic system with
a fixed parameter $\varepsilon >0\quad $

\begin{eqnarray}  \label{AEDFV}
\partial _{t}u^{\varepsilon }+\mathrm{div}\big(\mathbf{v}^{\varepsilon
}\;g(u^{\varepsilon })\big) &=&\varepsilon \;\Delta u^{\varepsilon }\text{%
\qquad in}\;\Omega _{T},  \label{AETU} \\
\tau \,\partial _{t}\mathbf{v}^{\varepsilon }-\nu \Delta \mathbf{v}%
^{\varepsilon }+h(u^{\varepsilon })\mathbf{v}^{\varepsilon } &=&-\nabla
p^{\varepsilon },\quad \mathrm{div}\left( \mathbf{v}^{\varepsilon }\right) =0%
\text{\quad in}\;\Omega _{T}  \label{AEPEV}
\end{eqnarray}%
jointly with the boundary-initial conditions%
\begin{eqnarray}
\varepsilon \frac{\partial u^{\varepsilon }}{\partial \mathbf{n}}%
+M(u^{\varepsilon }-u_{b}^{\varepsilon }) &=&0\text{\qquad and\qquad }%
\mathbf{v}^{\varepsilon }=\mathbf{b}\text{\quad on}\quad \Gamma _{T}  \notag
\\
(u^{\varepsilon },\mathbf{v}^{\varepsilon })|_{t=0} &=&(u_{0}^{\varepsilon },%
\mathbf{v}_{0})\quad \text{in }\Omega ,  \label{BC2}
\end{eqnarray}%
where $u_{b}^{\varepsilon },$ $u_{0}^{\varepsilon }$ are regularized
boundary-initial data satisfying suitable compatibility conditions. We
remark that $u_{b}^{\varepsilon }$, $u_{0}^{\varepsilon }$ converge strongly
in $L_{\,\mathrm{loc}}^{1}(\Gamma _{T})$ and $L_{\,\mathrm{loc}}^{1}(\Omega
),$ respectively, to $u_{b}$, $u_{0}.$

In the section \ref{SSKL2} we establish the following result.

\begin{proposition}
\label{PPA} For each $\varepsilon >0$, there exists a unique solution $%
(u^{\varepsilon },\mathbf{v}^{\varepsilon })$ of the system \eqref{AETU}--%
\eqref{BC2}, which has the following regularity $u^{\varepsilon }$ $\in
L^{\infty }(0,T;H^{1}(\Omega ))\;\cap \;L^{2}(0,T;H^{2}(\Omega ))$ and $%
\mathbf{v}^{\varepsilon }\in L^{2}(0,T;\mathbf{V}^{1}(\Omega ))\;\cap
\;H^{1}(0,T;\mathbf{V}^{-1}(\Omega ))$ satisfying
\begin{equation}
\varepsilon \Vert \nabla u^{\varepsilon }\Vert _{L^{2}(\Omega
_{T})}^{2}\leqslant C\qquad \text{and}\qquad 0\leqslant u^{\varepsilon
}\leqslant 1\quad \text{a.e. on $\Omega _{T}$},  \label{EUU}
\end{equation}%
\begin{equation}
\Vert \sqrt{\tau }\mathbf{v}^{\varepsilon }\Vert _{C([0,T];\mathbf{V}%
^{0}(\Omega ))}+\Vert \mathbf{v}^{\varepsilon }\Vert _{L^{2}(0,T;\mathbf{V}%
^{1}(\Omega ))}+\Vert \tau \mathbf{v}^{\varepsilon }\Vert _{H^{1}(0,T;%
\mathbf{V}^{-1}(\Omega ))}\leqslant C,  \label{EUU2}
\end{equation}%
where $C$ is a positive constant independent of $\varepsilon $ (and $\tau ).$
\end{proposition}

\bigskip

\begin{remark}
\label{rem2} After obtaining that \textbf{\ }$0\leqslant u^{\varepsilon
}\leqslant 1$ we can consider that $g(s):=g(0),$ $h(s):=h(0)$ for any $s<0$
\ and \ $g(s):=g(1),$ $h(s):=h(1)$ for any $s>1.$
\end{remark}

\subsection{\protect\bigskip The limit transition on $\protect\varepsilon %
\rightarrow 0$}

\bigskip

\label{SSKL}

\bigskip

In this section we are concerned to pass to the limit in \eqref{AETU}--%
\eqref{AEPEV} as $\varepsilon \rightarrow 0$. Since this problem is
non-linear on $u^{\varepsilon }$, the estimates \eqref{EUU}-\eqref{EUU2} are
not sufficient to take the limit transition on $\varepsilon $ as it goes to $%
0$. In fact, we need a strong convergence of a subsequence for the family $%
\{u^{\varepsilon }\}_{\varepsilon >0}$. Then, to derive the strong
convergence for $u^{\varepsilon }$, we use the Theory of Kinetic Formulation
as introduced by Lions, Perthame and Tadmor \cite{LPT1}-\cite{LPT2}, \cite%
{perthame2}. Here, we are going to follow closer Perthame and Dalibard \cite%
{PD}. That is, first we take the Kinetic formulation of \eqref{AETU}--%
\eqref{AEPEV}, then we pass to the weak limit. Finally, the information that
the initial-boundary conditions converge strongly, we are able to show the
strong convergence of $u^{\varepsilon }$.

\medskip

\subsubsection{\protect\vspace{1pt}The main idea of the limit transition.
Sketch of the proof}

\vspace{1pt}

Let $(\eta (u),q(u))$ be an entropy pair for \eqref{ETU}. Then, we have in
distribution sense
\begin{equation*}
\partial _{t}\eta (u^{\varepsilon })+\mathrm{div}(\mathbf{v}^{\varepsilon
}q(u^{\varepsilon }))-\varepsilon \;\Delta \eta (u^{\varepsilon
})=-\varepsilon \;\eta ^{\prime \prime }(u)\;|\nabla \eta (u)|^{2}\leqslant
0,
\end{equation*}%
since $\eta $ is a convex function. For instance, we could take the entropy
pair $(\eta (u),q(u))=$ $\mathbf{F}^{+}(u,v)$\ \ for all $v\in \mathbb{R},$
defined by \eqref{KE1}. Then, we have in sense of distributions
\begin{equation}
\partial _{t}|u^{\varepsilon }-v|^{+}+\mathrm{div}\left[ \mathbf{v}%
^{\varepsilon }\,\mathrm{sgn}^{+}\left( u^{\varepsilon }-v\right)
(g(u^{\varepsilon })-g(v))\right] -\varepsilon \;\Delta |u^{\varepsilon
}-v|^{+}=-m^{\varepsilon },  \label{ETK}
\end{equation}%
where $m^{\varepsilon }$ is a real nonnegative Radon measure$.$

\medskip If we differentiate in the distribution sense \eqref{ETK} with
respect to $v$, we get (as now a standard procedure in the kinetic theory)
the following transport equation
\begin{equation}
\partial _{t}f^{\varepsilon }+g^{\prime }(v)\,\mathbf{v}^{\varepsilon }\cdot
\nabla f^{\varepsilon }-\varepsilon \;\Delta f^{\varepsilon }=\partial
_{v}m^{\varepsilon },  \label{EDTK}
\end{equation}%
where $f^{\varepsilon }(t,\mathbf{x},v):=\mathrm{sgn}^{+}\left(
u^{\varepsilon }(t,\mathbf{x})-v\right) .$ \ Let us point out that
\begin{equation*}
0\leqslant f^{\varepsilon }(t,\mathbf{x},v)\leqslant 1\qquad \text{in $%
\Omega _{T}\times \mathbb{R}.$}
\end{equation*}%
Later on we show that $m^{\varepsilon }$ is uniformly bounded with respect
to $\varepsilon $, hence using \eqref{EUU}-\eqref{EUU2} there exist
subsequences of the families $\{m^{\varepsilon },f^{\varepsilon },\mathbf{v}%
^{\varepsilon }\}$ and a real nonnegative Radon measure $m=m(t,\mathbf{x}%
,v), $ functions $f\in L^{\infty }(\Omega _{T}\times \mathbb{R})$ and $%
\mathbf{v}\in L^{2}(0,T;\mathbf{V}^{1}(\Omega ))$, such that
\begin{align*}
m^{\varepsilon }& \rightarrow m\qquad \text{weakly in ${\mathcal{M}}(\Omega
_{T}\times \mathbb{R})$}, \\
f^{\varepsilon }& \rightarrow f\qquad \text{$\star $-weakly in $L^{\infty
}(\Omega _{T}\times \mathbb{R})$}, \\
\mathbf{v}^{\varepsilon }& \rightarrow \mathbf{v}\qquad \text{strongly in $%
L^{2}(\Omega _{T})$}.
\end{align*}%
Since \eqref{EDTK} is linear, it follows that
\begin{equation}
\partial _{t}f+g^{\prime }(v)\;\mathbf{v}\cdot \nabla f=\partial _{v}m\qquad
\text{in ${\mathcal{D}}^{\prime }(\Omega _{T}\times \mathbb{R})$}.
\label{EDTKL}
\end{equation}%
Accounting the initial boundary conditions for $f^{\varepsilon },$ we also
obtain
\begin{equation}
f=\mathrm{sgn}^{+}\left( u_{0}-v\right) \qquad \text{for $t=0$}\qquad \text{%
and}\qquad f=\mathrm{sgn}^{+}\left( u_{b}-v\right)  \label{bc}
\end{equation}%
on the influx part of $\Gamma _{T}\times \mathbb{R},$ \ i.e. where $%
g^{\prime }(v)\mathbf{b}_{\mathbf{n}}<0.$\ \ By the regularity of the
velocity field $\mathbf{v}\in L^{2}(0,T;\mathbf{V}^{1}(\Omega ))$, we can
use the theory for transport equations, introduced by DiPerna-Lions \cite%
{DiP-L}, and deduce that the solution of \ \eqref{EDTKL}-\eqref{bc} takes
values equals only to $0$ and $1\quad $on $\Omega _{T}\times \mathbb{R}$.
Since $f(\cdot ,\mathbf{\cdot },v)$ is a monotone function on $v$ (as a
limit of $f^{\varepsilon }(\cdot ,\mathbf{\cdot },v)$ \ being monotone one
too), we have
\begin{equation*}
f=\mathrm{sgn}^{+}\left( z(t,\mathbf{x})-v\right) \qquad \text{for some \ }%
z=z(t,\mathbf{x}).
\end{equation*}%
Finally, simplest considerations will apply that $z(t,\mathbf{x})\equiv u(t,%
\mathbf{x})$ and we have a strong convergence of \ $u^{\varepsilon }$ to $u,$
that ends the proof of our convergence result.

\vspace{1pt}

\vspace{1pt}

\subsubsection{Proof of Proposition\textbf{\ }\protect\ref{PPA}}

\label{SSKL2}

\vspace{1pt}Let $(\eta (u),q(u))$ be an entropy pair, satisfying the
condition
\begin{equation}
|q(u)|\leqslant K\eta (u)\quad \text{for}\quad u\in \mathbb{R}.
\label{ineqeta}
\end{equation}%
Both the pairs $\mathbf{F}^{\pm }(u,v)$\ \ for any $v\in \mathbb{R},$
defined by \eqref{KE1}, as the pair $\eta (u)=u^{2}$, $q(u)=\int_{0}^{u}2s\
g^{\prime }(s)\ ds$ fulfill this condition.

If we multiply \eqref{AETU} by $\eta ^{\prime }(u^{\varepsilon })\phi $ with
a function $\phi \in C_{0}^{\infty }((-\infty ,T)\times \mathbb{R}^{d})$ and
integrate on $\Omega _{T}$, we obtain
\begin{align}
\iint_{\Omega _{T}}& \left[ \eta (u^{\varepsilon })\phi
_{t}+q(u^{\varepsilon })\left( \mathbf{v}^{\varepsilon }\cdot \nabla \right)
\phi -\varepsilon \,\left( \nabla \phi \cdot \nabla \eta (u^{\varepsilon
})\right) \right] \;d\mathbf{x}dt  \notag \\
& +\int_{\Omega }\eta (u_{0}^{\varepsilon })\,\phi (0)\,d\mathbf{x}%
+\int_{\Gamma _{T}}M\eta (u_{b}^{\varepsilon })\,\phi \;d\mathbf{x}%
dt=m_{\eta }^{\varepsilon }(\phi ),  \label{EPB}
\end{align}%
where
\begin{align}
m_{\eta }^{\varepsilon }(\phi ):=& \iint_{\Omega _{T}}\varepsilon \eta
^{\prime \prime }|\nabla u^{\varepsilon }|^{2}\phi \ dtd\mathbf{x}%
+\int_{\Gamma _{T}}\{\mathbf{b_{\mathbf{n}}} \; q(u^{\varepsilon })+M
\eta(u^{\varepsilon })  \notag \\
& +\frac{1}{2}M\eta ^{\prime \prime }(r)(u_{b}^{\varepsilon }-u^{\varepsilon
})^{2}\}\phi \ dtd\mathbf{x}.  \label{m}
\end{align}%
Here we used that $\eta (u_{b}^{\varepsilon })=\eta (u^{\varepsilon })+\eta
^{\prime }(u^{\varepsilon })(u_{b}^{\varepsilon }-u^{\varepsilon })+\frac{%
\eta ^{\prime \prime }(r)}{2}(u_{b}^{\varepsilon }-u^{\varepsilon })^{2}$
for some function $r$ with values between $u^{\varepsilon }$ and $%
u_{b}^{\varepsilon }$ a.e. on $\Gamma _{T}.$ Let us observe that
\begin{equation}
m_{\eta }^{\varepsilon }(\phi )\geqslant \int_{\Omega _{T}}\varepsilon \eta
^{\prime \prime }|\nabla u^{\varepsilon }|^{2}\phi \ dtd\mathbf{x}\geqslant
0\qquad \text{for \quad }\phi \geqslant 0.  \label{m1}
\end{equation}%
\vspace{10pt}

Choosing in \eqref{EPB} $\phi (t,\mathbf{x}):=1-\zeta {\ _{\delta }}%
(t-t_{0}) $ for $t_{0}\in (0,T)$ with%
\begin{equation}
\zeta _{\delta }(s):=\left\{
\begin{array}{l}
0,\hspace{20pt}\text{if }s<0\text{\quad and\quad }1,\hspace{20pt}\text{if }%
s>\delta ,\text{ } \\
\frac{s}{\delta },\hspace{20pt}\text{if }0\leqslant s\leqslant \delta%
\end{array}%
\right.  \label{kzi}
\end{equation}%
and passing to the limit on $\delta \rightarrow 0^{+}$, we derive%
\begin{equation*}
\int_{\Omega }\eta (u^{\varepsilon }){\ }d\mathbf{x}(t_{0})+\int_{0}^{t_{0}}%
\int_{\Omega }\varepsilon \eta ^{\prime \prime }|\nabla u^{\varepsilon }|^{2}%
{\ }dtd\mathbf{x\leqslant }\int_{\Omega }\eta (u_{0}^{\varepsilon })\,\,d%
\mathbf{x}+\int_{0}^{t_{0}}\int_{\Gamma }M\eta (u_{b}^{\varepsilon })\,\,dtd%
\mathbf{x}.
\end{equation*}%
Hence taking $\eta (u)=|u|^{-}$\ ( $\eta =|u-1|^{+}$ and $\eta =u^{2},$\
consistently) \ in this inequality, we obtain the estimates (\ref{EUU}) by
the regularity assumptions \eqref{regular1}. The regularity $u^{\varepsilon
} $ $\in L^{\infty }(0,T;H^{1}(\Omega ))\;\cap \;L^{2}(0,T;H^{2}(\Omega ))$
follows from the well-known theory for parabolic type equations (see
Ladyzhenskaya et all \cite{LSU68}).

\bigskip

Now let us consider the quasi-stationary Stokes type equations
\begin{equation}
\begin{cases}
-\nu \Delta \mathbf{v}_{b}=-\nabla p_{b},\qquad \mathrm{div}(\mathbf{v}%
_{b})=0\quad \text{in }\Omega _{T}, \\
\mathbf{v}_{b}=\mathbf{b}\quad \text{on }\Gamma _{T}.%
\end{cases}
\label{v1}
\end{equation}%
In view of Proposition \ref{prop2} and the assumption \eqref{regular2}, the
solution $\mathbf{v}_{b}$ of this problem exists and fulfills the estimate%
\begin{equation}
\Vert \mathbf{v}_{b}\Vert _{L^{2}(0,T;\mathbf{H}^{1}(\Omega ))}+\Vert
\partial _{t}\mathbf{v}_{b}\Vert _{\mathbf{L}^{2}(\Omega _{T})}\leqslant C.
\label{v2}
\end{equation}%
Therefore taking the difference between \eqref{AEPEV} and \eqref{v1}, we
have that the function $\mathbf{w}^{\varepsilon }=\mathbf{v}^{\varepsilon }-%
\mathbf{v}_{b}$ satisfies
\begin{eqnarray*}
\tau \,\partial _{t}\mathbf{w}^{\varepsilon }-\nu \Delta \mathbf{w}%
^{\varepsilon }+h(u^{\varepsilon })\mathbf{w}^{\varepsilon } &=&-\nabla
(P^{\varepsilon })+\mathbf{f}^{\varepsilon },\qquad \mathrm{div}\left(
\mathbf{w}^{\varepsilon }\right) =0, \\
\mathbf{w}^{\varepsilon }\big|_{\Gamma _{T}} &=&0,\qquad \mathbf{w}%
^{\varepsilon }\big|_{t=0}=\mathbf{v}_{0}-\mathbf{v}_{b}\big|_{t=0},
\end{eqnarray*}%
with $P^{\varepsilon }:=p^{\varepsilon }-p_{b},$ $\ \mathbf{f}^{\varepsilon }%
\mathbf{:=-}\tau \partial _{t}\mathbf{v}_{b}-h(u^{\varepsilon })\mathbf{v}%
_{b}.$\ \ Let us point that the solvability of the above system can been
shown as in \cite{Lad69}, \cite{Temam}. If we multiply the first equation in
this system by $\mathbf{w}^{\varepsilon }$ and integrate over $\Omega ,$ we
obtain
\begin{align*}
\frac{d}{dt}\left( \frac{\tau }{2}\Vert \mathbf{w}^{\varepsilon }\Vert _{%
\mathbf{L}^{2}(\Omega ))}^{2}\right) +\nu \Vert \nabla \mathbf{w}%
^{\varepsilon }\Vert _{\mathbf{L}^{2}(\Omega )}^{2}& \leqslant \int_{\Omega
}|(\mathbf{f}^{\varepsilon }\cdot \mathbf{w}^{\varepsilon })|\,d\mathbf{x} \\
& \leqslant \frac{\nu }{2}\Vert \nabla \mathbf{w}^{\varepsilon }\Vert _{%
\mathbf{L}^{2}(\Omega )}^{2}+C\Vert \mathbf{f}^{\varepsilon }\Vert _{\mathbf{%
L}^{2}(\Omega )}^{2},
\end{align*}%
where we have used Poincaré's inequality.

Then, using (\ref{regular2}), (\ref{v2}) we deduce
\begin{equation*}
\tau \Vert \mathbf{w}^{\varepsilon }\Vert _{L^{\infty }(0,T;\mathbf{L}%
^{2}(\Omega ))}^{2}+\Vert \mathbf{w}^{\varepsilon }\Vert _{L^{2}(0,T;\mathbf{%
V}^{1}(\Omega ))}^{2}\leqslant C
\end{equation*}%
with some constant $C$ independent of $\varepsilon $ (and $\tau ).$ Hence in
view of the weak formulation (\ref{V}) of \eqref{AEPEV} and Lemmas 1.2-1.4,
p.176 of \cite{Temam}, we get that $\mathbf{w}^{\varepsilon }\in C([0,T];%
\mathbf{V}^{0}(\Omega ))$ and the estimate (\ref{EUU2}).

Finally, with the help of derived estimates
(\ref{EUU})-(\ref{EUU2}), we can apply Leray-Schauder's fixed
point argument (as now a standard procedure) and get the
solvability of the approximated system \eqref{AETU}--\eqref{BC2}.
\qquad $\blacksquare $

\vspace{1pt}

\vspace{1pt}

\subsubsection{\protect\vspace{1pt}Rigorous proof of the limit transition}

\vspace{1pt}\label{SD}

Now, if we take in \eqref{EPB} the entropy pair $\mathbf{F}^{+}(u,v)$\ \ for
all $v\in \mathbb{R},$ then we see that the function $f^{\varepsilon }(t,%
\mathbf{x},v)=\,\mathrm{sgn}^{+}\left( u^{\varepsilon }-v\right) $ satisfies
for all nonnegative function $\phi \in C_{0}^{\infty }((-\infty ,T)\times
\mathbb{R}^{d})$, the following identity
\begin{align}
\iint_{\Omega _{T}}& \left\{ \int_{v}^{1}f^{\varepsilon }(t,\mathbf{x},s)%
\left[ \phi _{t}+g^{\prime }(s)\left( \mathbf{v}^{\varepsilon }\cdot \nabla
\right) \phi \right] \;ds-\,\varepsilon \nabla \phi \cdot \nabla
\,|u^{\varepsilon }-v|^{+}\right\} \;d\mathbf{x}dt  \notag \\
& +\int_{\Omega }|u_{0}^{\varepsilon }-v|^{+}\ \phi (0)\,d\mathbf{x+}%
\int_{\Gamma _{T}}M|u_{b}^{\varepsilon }-v|^{+}\;\phi \;d\mathbf{x}%
dt=m_{+}^{\varepsilon }(\phi )\geqslant 0,  \label{EPB+}
\end{align}%
where $m_{+}^{\varepsilon }:=m_{|u^{\varepsilon }-v|^{+}}^{\varepsilon }$ \
(see \eqref{m} and \eqref{m1}).

\bigskip Further, we have for any $G\in C^{1}([0,1])$, with $G(0)=0$ that%
\begin{eqnarray}
G(u^{\varepsilon }) &=&\int_{0}^{1}G^{\prime }(s)f^{\varepsilon }(\cdot ,%
\mathbf{\cdot },s)\;ds\qquad \text{ a. e. in $\Omega _{T},$}  \notag \\
0 &\leqslant &f^{\varepsilon }\leqslant 1\qquad \text{on $\Omega _{T}\times
\mathbb{R}$},\qquad \text{$f^{\varepsilon }(t,\mathbf{x},v)=\left\{ %
\begin{aligned} 1, &\quad \text{for $v \leq 0$}, \\ 0, &\quad \text{for $v
\geq 1$}, \end{aligned}\right. $}  \notag \\
\partial _{v}f^{\varepsilon } &\leqslant &0\qquad \text{in ${\mathcal{D}}%
^{\prime }(\Omega _{T}\times \mathbb{R})$}.  \label{u-f}
\end{eqnarray}%
Let us choose in \eqref{EPB+} \ $\phi :=1-\zeta _{\delta }(t-T+\delta )$ and
then, passing to the limit as $\delta \rightarrow 0^{+},$ we get%
\begin{equation}
m_{+}^{\varepsilon }(1)\leqslant \int_{\Omega }|u_{0}^{\varepsilon
}-v|^{+}\;d\mathbf{x+}\int_{\Gamma _{T}}M|u_{b}-v|^{+}\;d\mathbf{x}dt,
\label{m-estimate}
\end{equation}%
hence the Riesz representation theorem implies that the real positive Radon
measure $m_{+}^{\varepsilon }$ is well defined on $\overline{\Omega }%
_{T}\times \mathbb{R},$ such that
\begin{align}
& m_{+}^{\varepsilon }(\cdot ,\cdot ,v)=0\qquad \text{for any }v>1 \quad
\text{on }\Omega _{T},\; \text{and }  \notag \\
& \iint_{\overline{\Omega }_{T}}m_{+}^{\varepsilon }(t,\mathbf{x},v)\;d%
\mathbf{x}dt\leqslant C(v)\quad \text{ for all finite }v,  \label{m-property}
\end{align}%
where $C(v)$ is a positive constant independent of $\varepsilon ,$ but could
depend on $v$. By a similar way as \eqref{EPB+} has been derived, if we take
the entropy pair $\mathbf{F}^{-}(u,v)$\ \ for all $v\in \mathbb{R},$ defined
by \eqref{KE1}, we obtain
\begin{align}
\iint_{\Omega _{T}}& \left\{ \int_{0}^{v}(1-f^{\varepsilon }(t,\mathbf{x},s))%
\left[ \phi _{t}+g^{\prime }(s)\left( \mathbf{v}^{\varepsilon }\cdot \nabla
\right) \phi \right] \;dv-\,\varepsilon \nabla \phi \cdot \nabla
\,|u^{\varepsilon }-v|^{-}\right\} \;d\mathbf{x}dt  \notag \\
& +\int_{\Omega }|u_{0}^{\varepsilon }-v|^{-}\;\phi (0)\,\;d\mathbf{x+}%
\int_{\Gamma _{T}}M|u_{b}^{\varepsilon }-v|^{-}\;\phi \;d\mathbf{x}%
dt=m_{-}^{\varepsilon }(\phi )\geqslant 0,  \label{EPB-}
\end{align}%
where $m_{-}^{\varepsilon }:=m_{|u^{\varepsilon }-v|^{-}}^{\varepsilon }.$ \
Moreover we have that the real positive Radon measure $m_{-}^{\varepsilon },$
defined on $\overline{\Omega }_{T}\times \mathbb{R}$,\ fulfills the
following properties
\begin{align}
m_{-}^{\varepsilon }(1)& \leqslant \int_{\Omega }|u_{0}^{\varepsilon
}-v|^{-}\;d\mathbf{x+}\int_{\Gamma _{T}}M|u_{b}^{\varepsilon }-v|^{-}\;d%
\mathbf{x}dt,  \notag \\
m_{-}^{\varepsilon }(\cdot ,\cdot ,v)& =0\qquad \text{for any }v<0\text{ and
on }\overline{\Omega }_{T},  \notag \\
\iint_{\overline{\Omega }_{T}}m_{-}^{\varepsilon }(t,\mathbf{x},v)\;d\mathbf{%
x}dt& \leqslant C(v)\quad \text{ for all finite }v.  \label{m-property2}
\end{align}

In view of Proposition \ref{PPA} and \eqref{m-property}, \eqref{m-property2}%
, \ there exist subsequences of $f^{\varepsilon }$, $\mathbf{v}^{\varepsilon
}$, $m^{\varepsilon }$ and the functions
\begin{equation}
f\in L^{\infty }(\Omega _{T}\times \mathbb{R}),\quad \mathbf{v}\in L^{2}(0,T;%
\mathbf{V}^{1}(\Omega )),  \label{reg}
\end{equation}%
and a real nonnegative Radon measure $m=m(t,\mathbf{x},v)$, such that
\begin{align*}
f^{\varepsilon }& \rightarrow f\qquad \text{$\star $-weakly in $L^{\infty
}(\Omega _{T}\times \mathbb{R})$},\qquad \\
\mathbf{v}^{\varepsilon }& \rightarrow \mathbf{v},\qquad \varepsilon
\,\nabla u^{\varepsilon }\rightarrow 0\qquad \text{strongly in $L^{2}(\Omega
_{T})$}, \\
m_{\pm }^{\varepsilon }& \rightarrow m_{\pm }\qquad \text{weakly in ${%
\mathcal{M}}_{loc}^{+}(\overline{\Omega }_{T}\times \mathbb{R}).$}
\end{align*}%
Now, for any nonnegative function $\phi \in C_{0}^{\infty }((-\infty
,T)\times \mathbb{R}^{d+1})$ the following integral inequalities fulfills
\begin{align}
& \iint_{\Omega _{T}}\int_{v}^{1}f(t,\mathbf{x},s)\left[ \phi _{t}+g^{\prime
}(s)\left( \mathbf{v}\cdot \nabla \right) \phi \right] \;dsd\mathbf{x}%
dt+\int_{\Omega }|u_{0}-v|^{+}\phi (0)\;d\mathbf{x}  \notag \\
& +\int_{\Gamma _{T}}M\ |u_{b}-v|^{+}\phi \;d\mathbf{x}dt=m_{+}(\phi
)=:\int_{\Omega _{T}}m_{+}(t,\mathbf{x},v)\phi \;d\mathbf{x}dt\geqslant 0,
\label{EPB+L}
\end{align}%
and
\begin{align}
& \iint_{\Omega _{T}}\int_{0}^{v}(1-f(t,\mathbf{x},s))\left[ \phi
_{t}+g^{\prime }(s)\left( \mathbf{v}\cdot \nabla \right) \phi \right] \;dsd%
\mathbf{x}dt+\int_{\Omega }|u_{0}-v|^{-}\phi (0)\;d\mathbf{x}  \notag \\
& +\int_{\Gamma _{T}}M\ |u_{b}-v|^{-}\phi \;d\mathbf{x}dt=m_{-}(\phi
)=:\int_{\Omega _{T}}m_{-}(t,\mathbf{x},v)\phi \;d\mathbf{x}dt\geqslant 0.
\label{EPB-L}
\end{align}%
Moreover, we have for any $G\in C^{1}([0,1])$, with $G(0)=0$ that%
\begin{eqnarray}
G(u) &=&\int_{0}^{1}G^{\prime }(s)f(\cdot ,\mathbf{\cdot },s)\;ds\qquad
\text{ a.e. in $\Omega _{T},$}  \notag \\
0 &\leqslant &f\leqslant 1\qquad \text{on $\Omega _{T}\times \mathbb{R}$}%
,\qquad \text{$f(t,\mathbf{x},v)=\left\{ \begin{aligned} 1, &\quad \text{for
$v \leq 0$}, \\ 0, &\quad \text{for $v \geq 1$}, \end{aligned}\right. $}
\notag \\
\partial _{v}f &\leqslant &0\qquad \text{in ${\mathcal{D}}^{\prime }(\Omega
_{T}\times \mathbb{R}),$}  \label{for u}
\end{eqnarray}%
and
\begin{align}
\iint_{\overline{\Omega }_{T}}m_{\pm }(t,\mathbf{x},v)\;d\mathbf{x}dt&
\leqslant C(v)\quad \text{ for all finite }v,  \notag \\
m_{+}(\cdot ,\cdot ,v)& =0\qquad \text{for any }v>1\text{ and on }\overline{%
\Omega }_{T},  \notag \\
m_{-}(\cdot ,\cdot ,v)& =0\qquad \text{for any }v<0\text{ and on }\overline{%
\Omega }_{T},  \notag \\
m_{\pm }(\cdot ,\cdot ,v)& \in C(\mathbb{R};\text{${\mathcal{M}}^{+}(%
\overline{\Omega }_{T}\times \mathbb{R})$}),  \label{for m}
\end{align}%
the continuity of $m_{\pm }(\cdot ,\cdot ,v)$ on $v$ \ follows from the left
parts of \eqref{EPB+L}, \eqref{EPB-L}.

\medskip Finally, taking in \eqref{EPB+L} and \eqref{EPB-L} $\phi =\partial
_{v}\psi $, with $\psi $ being a nonnegative function in $C_{0}^{\infty
}(\Omega _{T}\times \mathbb{R})$, integrating by parts on $v$, we obtain
that $f$ satisfies the following transport equations%
\begin{eqnarray}
\iint_{\Omega _{T}\times \mathbb{R}}f\left[ \psi _{t}+g^{\prime }(v)\left(
\mathbf{v}\cdot \nabla \right) \psi \right] -m_{+}\,\partial _{v}\psi \;dvd%
\mathbf{x}dt &=&0,  \notag \\
\iint_{\Omega _{T}\times \mathbb{R}}(1-f)\left[ \psi _{t}+g^{\prime
}(v)\left( \mathbf{v}\cdot \nabla \right) \psi \right] +m_{-}\,\partial
_{v}\psi \;dvd\mathbf{x}dt &=&0,  \label{f}
\end{eqnarray}%
respectively.

\bigskip Now, let us study the trace concept on the initial-boundary terms.

\begin{proposition}
\label{TRACES} The function $f=f(t,\mathbf{x},v)$ has the trace $f^{0}$ at
the time $t=0,$ such that
\begin{equation*}
f^{0}=f^{0}(\mathbf{x},v)\equiv \lim_{\delta \rightarrow 0^{+}}\frac{{1}}{%
\delta }{\int }_{\!\!\!0}^{\delta }f(s,\mathbf{x},v)\;ds
\end{equation*}%
and
\begin{equation}
f^{0}=\left( f^{0}\right) ^{2}\qquad \text{a. a. on $\Omega \times \mathbb{R}
$}.  \label{TIC}
\end{equation}%
The function $f=f(t,\mathbf{x},v)$ has the trace $f^{b}$ on $\Gamma
_{T}\times \mathbb{R}$, such that
\begin{equation*}
f^{b}=f^{b}(t,\mathbf{x},v)\equiv \lim_{\delta \rightarrow 0^{+}}{\frac{{1}}{%
\delta }\int }_{\!\!\!0}^{\delta }f(t,\mathbf{x}-s\,n(\mathbf{x}),v)\;ds
\end{equation*}%
for a. a. $(t,\mathbf{x},v)\in \Gamma _{T}\times \mathbb{R},$ where $%
\;g^{\prime }(v)\mathbf{b}_{\mathbf{n}}(t,\mathbf{x})\not=0\;$ and
\begin{equation}
f^{b}=\left( f^{b}\right) ^{2}  \label{2.131}
\end{equation}%
for a. a. $(t,\mathbf{x},v)\in \Gamma _{T}\times \mathbb{R},$ where $%
\;g^{\prime }(v)\mathbf{b}_{\mathbf{n}}(t,\mathbf{x})<0.$
\end{proposition}

\textbf{Proof.} First, let $\gamma \in C_{0}^{\infty }(\mathbb{R})$ be a
fixed function. Then by \eqref{f} the vector function
\begin{equation*}
\mathbf{f}_{\gamma }:=\Big(\int_{R}f(\cdot ,\cdot ,v)\gamma (v)dv,\quad
\mathbf{v}\int_{R}g^{\prime }(v)f(\cdot ,\cdot ,v)\gamma (v)dv\Big)\in
\mathbf{L}^{2}(\Omega _{T}),
\end{equation*}%
and it follows that
\begin{equation*}
\mathrm{div}_{(t,\mathbf{x})}(\mathbf{f}_{\gamma })=-\int_{R}m_{+}(\cdot
,\cdot ,v)\gamma ^{\prime }(v)\;dv\in {\mathcal{M}}(\Omega _{T})
\end{equation*}%
having a finite total variation $|\mathrm{div}_{(t,\mathbf{x})}f_{\gamma
}|(\Omega _{T})<\infty $, in view of \eqref{for m}. Let $\Sigma $ and $%
\mathbf{n}_{(t,\mathbf{x})}$\ be the boundary of $\overline{\Omega }_{T}$
and the unitary normal to $\Sigma ,$ respectively. By Theorem 2.1 in Chen \&
Frid \cite{CF}, $\mathbf{f}_{\gamma }\cdot \mathbf{n}_{(t,\mathbf{x})}$ is a
continuous linear functional over $H^{1/2}(\Sigma )\cap L^{\infty }(\Sigma
). $

1.1 Now, if we take $\varphi \in C_{0}^{\infty }((-\infty ,T)\times \Omega )$%
, then
\begin{align*}
<\mathbf{f}_{\gamma }\cdot \mathbf{n}_{(t,\mathbf{x})}|_{t=0},\varphi >&
=\lim_{\delta \rightarrow 0^{+}}\frac{{1}}{\delta }{\int }_{\!\!\!0}^{\delta
}\int_{\Omega }\left[ \int_{\mathbb{R}}f(t,\mathbf{x},v)\gamma (v)dv\right]
\varphi (t,\mathbf{x})\;d\mathbf{x}dt \\
& =\lim_{\delta \rightarrow 0^{+}}\iint_{\Omega \times \mathbb{R}}\left[
\frac{{1}}{\delta }{\int }_{\!\!\!0}^{\delta }f(t,\mathbf{x},v)\;dt\right]
\gamma (v)\varphi (0,\mathbf{x})\;d\mathbf{x}dv.
\end{align*}%
Since $\ 0\leqslant \frac{{1}}{\delta }{\int }_{\!\!\!0}^{\delta }f(t,\cdot
,\cdot )\;dt\leqslant 1$ \quad on $\Omega \times \mathbb{R}$, using the
dominated convergence theorem we derive the existence of $\lim_{\delta
\rightarrow 0^{+}}\frac{{1}}{\delta }{\int }_{\!\!\!0}^{\delta }f(s,\cdot
,\cdot )\;ds$, which we denote by $f^{0}.$ It is obvious that $0\leqslant
f^{0}\leqslant 1$\quad\ on $\Omega \times \mathbb{R}$ and $\mathbf{f}%
_{\gamma }\cdot \mathbf{n}_{(t,\mathbf{x})}|_{t=0}=\mathbf{<f^{0},\gamma >}$
a.e. on $\mathbf{\Omega .}$ Since $\mathbf{\gamma }$ is an arbitrary
function,  therefore we can simply denote \textbf{$\mathbf{f}\cdot \mathbf{n}%
_{(t,\mathbf{x})}|_{t=0}=f^{0}.$}

\medskip 1.2. Let us take a nonnegative function $\psi \in C_{0}^{\infty
}(\Omega )$ and set $\varphi (t,\mathbf{x})=\zeta _{\delta }(t)\psi (\mathbf{%
x})$ in inequalities \eqref{EPB+L}, \eqref{EPB-L} with $\zeta _{\delta }$
given by \eqref{kzi}.\ Then, we obtain, after passing to the limit as $%
\delta \rightarrow 0^{+}$, respectively
\begin{equation*}
\int_{\Omega }\psi (\mathbf{x})\left[ -\int_{v}^{1}f^{0}(\mathbf{x},s)\;dsd%
\mathbf{x}+|u_{0}(\mathbf{x})-v|^{+}\right] \;d\mathbf{x}\geqslant 0,
\end{equation*}%
and
\begin{equation*}
\int_{\Omega }\psi (\mathbf{x})\left[ -\int_{0}^{v}(1-f^{0}(\mathbf{x}%
,s))\;\;dsd\mathbf{x}+|u_{0}(\mathbf{x})-v|^{-}\right] \;d\mathbf{x}%
\geqslant 0.
\end{equation*}%
Since $\psi $ is an arbitrary nonnegative function, then \ for a. a. $%
\mathbf{x}\in \Omega $\ \ the 1$^{st}$ \ inequality implies that $f^{0}(%
\mathbf{x},v)=0\quad \ $if$\ v>u_{0}(\mathbf{x})$ \ \ \ and \ \ the 2$^{nd}$
\ one gives $f^{0}(\mathbf{x},v)=1\quad $ if $\ v<u_{0}(\mathbf{x})$, i.e.
we show \eqref{TIC}.

\medskip 2.1. Let $d(\mathbf{x}):=\min_{\mathbf{y}\in \Gamma }|\mathbf{x}-%
\mathbf{y}|$ be the distance function from $\mathbf{x}\in \overline{\Omega }$
\ to $\ \Gamma .$ \ Denoting by $\mathbf{x}_{s}:=\mathbf{x-}s\mathbf{n}(%
\mathbf{x}),$ for any $\mathbf{x}\in \Gamma $ and $s\in (0,\delta ),$ and
applying again the result of Theorem 2.1 in \cite{CF}, we have that for any $%
\psi \in C_{0}^{\infty }((0,T)\times \mathbb{R}^{d})$
\begin{align*}
& <\mathbf{f}_{\gamma }\cdot \mathbf{n}|_{\Gamma _{T}},\psi > \\
& =\lim_{\delta \rightarrow 0^{+}}\frac{{1}}{\delta }{\int }%
_{\!\!\!0}^{\delta }\int_{\Gamma _{T}}\left[ \int_{\mathbb{R}}g^{\prime }(v)(%
\mathbf{v\cdot }(-\nabla d)\mathbf{)}(t,\mathbf{x}_{s})\ f(t,\mathbf{x}%
_{s},v)\;\gamma (v)\ dv\right] \;\psi (t,\mathbf{x}_{s})\;d\mathbf{x}dtds \\
& =\iint_{\Gamma _{T}\times \mathbb{R}}g^{\prime }(v)\;\mathbf{b}_{\mathbf{n}%
}(t,\mathbf{x})\;\psi (t,\mathbf{x})\;\gamma (v)\left[ \lim_{\delta
\rightarrow 0^{+}}\frac{{1}}{\delta }{\int }_{\!\!\!0}^{\delta }f(t,\mathbf{x%
}_{s},v)\;ds\right] \;d\mathbf{x}dtdv.
\end{align*}%
In the last equality we have used that $\mathbf{v}(t,\cdot )\in \mathbf{V}%
^{1}(\Omega )$ for a. a. $t\in \lbrack 0,T]$ with Theorems 6.5.3, 6.5.4 of
\cite{Kuf}$;$\ \ \ $-\nabla d\equiv \mathbf{n}$ on $\Gamma \ $ with $\Gamma
\in C^{2}$\ and also the dominated convergence theorem applied for the
bounded sequence%
\begin{equation*}
0\leqslant \frac{{1}}{\delta }{\int }_{\!\!\!0}^{\delta }f(t,\mathbf{x}%
_{s},v)\;ds\leqslant 1\qquad \text{for a. a. }(t,\mathbf{x},v)\in \Gamma
_{T}\times \mathbb{R},
\end{equation*}%
\ where $g^{\prime }(v)\mathbf{b}_{\mathbf{n}}\not=0.$ Denoting by $f^{b}:=$
$\lim_{\delta \rightarrow 0^{+}}\frac{{1}}{\delta }{\int }_{\!\!\!0}^{\delta
}f(t,\mathbf{x}_{s},v)\;ds,$ we have $\mathbf{f}_{\gamma }\cdot \mathbf{n}%
|_{\Gamma _{T}}\equiv g^{\prime }(v)\mathbf{b}_{\mathbf{n}}\mathbf{<}%
f^{b},\gamma >.$ Since $\gamma $ is an arbitrary function, hence we can
denote\ $\mathbf{f}\cdot \mathbf{n}|_{\Gamma _{T}}=g^{\prime }(v)\mathbf{b}_{%
\mathbf{n}}f^{b}.$  The function $f^{b}$ fulfills
\begin{equation}
0\leqslant f^{b}\leqslant 1\qquad \text{on $\Gamma _{T}\times \mathbb{R}$}%
\qquad \text{and}\qquad \partial _{v}f^{b}\leqslant 0\qquad \text{in ${%
\mathcal{D}}^{\prime }(\Gamma _{T}\times \mathbb{R})$}.  \label{TBC}
\end{equation}%
\

\vspace{1pt}

2.2. \ Let us take a positive function $\psi \in C_{0}^{\infty }((0,T)\times
\mathbb{R}^{d})$. If we set $\ \phi (t,\mathbf{x})=(1-\zeta _{\delta }(d(%
\mathbf{x})))\ \psi (t,\mathbf{x})\;$\ \ in \eqref{EPB+L}, \eqref{EPB-L}
with $\zeta _{\delta }$ defined in \eqref{kzi} \ and pass to the limit as $%
\delta \rightarrow 0^{+}$, we obtain respectively
\begin{equation*}
\int_{\Gamma _{T}}\psi \left[ \int_{v}^{1}g^{\prime }(s)\mathbf{b}_{\mathbf{n%
}}\;f^{b}(t,\mathbf{x},s)\ ds+M\ |u_{b}-v|^{+}\right] \;\;d\mathbf{x}%
dt\geqslant 0
\end{equation*}%
and
\begin{equation*}
\int_{\Gamma _{T}}\psi \left[ \int_{0}^{v}g^{\prime }(s)\mathbf{b}_{\mathbf{n%
}}(1-\;f^{b}(t,\mathbf{x},s))\;ds+M\ |u_{b}-v|^{-}\right] \;d\mathbf{x}%
dt\geqslant 0.
\end{equation*}%
Hence defining the functions
\begin{align*}
m_{+}^{b}& :=\int_{v}^{1}g^{\prime }(s)\mathbf{b}_{\mathbf{n}}\;f^{b}(t,%
\mathbf{x},s)\;ds+M\;|u_{b}-v|^{+}, \\
m_{-}^{b}& :=\int_{0}^{v}g^{\prime }(s)\mathbf{b}_{\mathbf{n}}\;\left(
1-f^{b}(t,\mathbf{x},s)\right) \;ds+M\;|u_{b}-v|^{-}
\end{align*}%
for $(t,\mathbf{x},v)\in \Gamma _{T}\times \mathbb{R}$, it is not difficult
to check that, the positive functions $m_{\pm }^{b}(t,\mathbf{x},v)\in
L^{2}(\Gamma _{T};W^{1,\infty }(\mathbb{R}))$ and satisfy
\begin{align}
g^{\prime }(v)\mathbf{b}_{\mathbf{n}}\;f^{b}& =-M\,\mathrm{sgn}^{+}\left(
u_{b}-v\right) -\partial _{v}m_{+}^{b}\quad \text{and}\quad  \notag \\
m_{+}^{b}& =0\quad \text{for $v\geqslant 1$};  \notag \\
g^{\prime }(v)\mathbf{b}_{\mathbf{n}}(1-f^{b})& =-M\,\mathrm{sgn}^{-}\left(
u_{b}-v\right) +\partial _{v}m_{-}^{b}\quad \text{and}\quad  \notag \\
m_{-}^{b}& =0\quad \text{for $v\leqslant 0.$}  \label{TBCM}
\end{align}

Due to (\ref{TBC})-(\ref{TBCM}) we have
\begin{align*}
0& \leqslant \int\limits_{\mathbb{R}}\left\vert g^{\prime }(v)\mathbf{b}_{%
\mathbf{n}}\right\vert
^{-}\,f^{b}(v)(1-f^{b}(v))\,dv=\int\limits_{0}^{u_{b}}f^{b}\{M\mathrm{sgn}%
^{-}\left( u_{b}-v\right) -\partial _{v}m_{-}^{b}\}\,dv \\
& +\int\limits_{u_{b}}^{1}\left\{ M\mathrm{sgn}^{+}\left( u_{b}-v\right)
+\partial _{v}m_{+}^{b}\right\} (1-f^{b})\,dv=-(f^{b}m_{-}^{b})|_{v=u_{b}-0}
\\
& +\int\limits_{0}^{u_{b}}\partial
_{v}f^{b}m_{-}^{b}\,dv-m_{+}^{b}(1-f^{b})|_{v=u_{b}+0}+\int%
\limits_{u_{b}}^{1}m_{+}^{b}\partial _{v}f^{b}dv\leqslant 0\text{\qquad a.e.
on }\Gamma _{T}.
\end{align*}%
A formal integration on $v$ by parts in the last identity can be justified
by mollifying the function $f^{b}$ and taking the limit transition on a
mollifying parameter. Therefore $f^{b}$\ satisfies \eqref{2.131}. \ $%
\blacksquare $

\bigskip

\begin{lemma}
\label{lemma33} We have $\ \ \ $%
\begin{equation}
f=f^{2}\qquad \text{a. e. in }\Omega _{T}\times \mathbb{R}.  \label{aA}
\end{equation}
\end{lemma}

\noindent \textbf{Proof.} \ The equations \eqref{f} are written as
\begin{eqnarray}
\partial _{t}f+\mathrm{div}_{\mathbf{x}}(g^{\prime }(v)\,\mathbf{v}\,f)
&=&\partial _{v}m_{+}  \notag \\
\partial _{t}(1-f)+\mathrm{div}_{\mathbf{x}}(g^{\prime }(v)\,\mathbf{v}%
\,(1-f)) &=&-\partial _{v}m_{-}\quad \text{ in \ }\mathcal{D}^{\prime
}(\Omega _{T}\mathbb{\times R}).  \label{eq2.9}
\end{eqnarray}%
We have that $\mathbf{v}\in L^{2}(0,T,\mathbf{V}^{1}(\Omega ))$ and $%
g^{\prime }(v)$\ is a constant with respect of the variables $(t,\mathbf{x}),
$\ such that $g^{\prime }\in L^{\infty }(\mathbb{R})$ by Remark \ref{rem2}.
Hence we can apply the renormalization theorem to the left parts of %
\eqref{eq2.9} \ \ (see for instance Theorem 4.3 in \cite{lellis}) and get,
that the function $F:=f(1-f)$ satisfies%
\begin{equation}
F_{t}+\mathrm{div}_{\mathbf{x}}(g^{\prime }(v)\mathbf{v}F)\leqslant 0\quad
\text{ in \ }\mathcal{D}^{\prime }(\mathbb{R\times }\Omega _{T}).
\label{eqq1}
\end{equation}%
It means that the equations in \eqref{eq2.9} are regularized on a parameter $%
\theta ,$ multiplied by $(1-f^{\theta })$ and $f^{\theta },$ respectively ($%
f^{\theta }$ being the regularization of $f)$\ and finally taken the limit
on $\theta \rightarrow 0.$ The inequality in \eqref{eqq1} follows from the
relation $\int\limits_{\mathbb{R}}\frac{\partial m_{+}^{\theta }}{\partial v}%
(1-f^{\theta })-\frac{\partial m_{-}^{\theta }}{\partial v}f^{\theta
}\,dv=\int\limits_{\mathbb{R}}(m_{+}^{\theta }+m_{-}^{\theta })\frac{%
\,\partial f^{\theta }}{\partial v}\,dv\leqslant 0$ \ in view of (\ref{for u}%
) and (\ref{for m}).

Now let us define%
\begin{equation*}
\psi :=(\zeta {_{\varepsilon }}(v+\varepsilon ^{-1})-\zeta {_{\varepsilon }}%
(v-\varepsilon ^{-1}))\psi _{1}^{\delta }(t)\psi _{2}^{\delta }(\mathbf{x})
\end{equation*}%
with $\psi _{1}^{\delta }(t):=(\zeta {_{\delta }}(t)-\zeta {_{\delta }}%
(t-t_{0}+\delta ))$ \ for \ $t_{0}\in (2\delta ,T)$ and $\psi _{2}^{\delta }(%
\mathbf{x}):=\zeta _{\delta }(d(\mathbf{x}))$\ for $\mathbf{x}\in \Omega .$
Choosing $\psi =\psi (v,t,\mathbf{x})$ as a test function in the respective
integral form of \eqref{eqq1} and taking the limit transition on $%
\varepsilon \rightarrow 0,$ with the help of (\ref{for u}) and (\ref{for m}%
), we get the inequality%
\begin{eqnarray}
&&\frac{1}{\delta }\int_{t_{0}-\delta }^{t_{0}}\int_{\Omega \times \mathbb{R}%
}F\ \psi _{2}^{\delta }(\mathbf{x})\ dtdvd\mathbf{x}\leqslant \frac{1}{%
\delta }\int_{0}^{\delta }\int_{\Omega \times \mathbb{R}}F\ \psi
_{2}^{\delta }(\mathbf{x})\ dtdvd\mathbf{x}  \notag \\
&&+\frac{1}{\delta }\int_{0\leqslant d(\mathbf{x})\leqslant \delta }\int_{{%
\mathbb{R}}\times \lbrack 0,T]}|g^{\prime }(v)(\mathbf{v\cdot }\nabla
d)|^{-}\ F\ \psi _{1}^{\delta }(t)\ dtdvd\mathbf{x}  \notag \\
&=&:C_{1}^{\delta }+C_{2}^{\delta }.  \label{ddF}
\end{eqnarray}

Due to the following simple inequality%
\begin{equation}
-\frac{1}{\delta }\int_{0}^{\delta }z^{2}(s)ds\leqslant -\left( \frac{1}{%
\delta }\int_{0}^{\delta }z(s)ds\right) ^{2},  \label{holder}
\end{equation}%
which is valid for any positive integrable function $z=z(s),$ \ we have that%
\begin{equation*}
C_{1}^{\delta }\leqslant \int_{\Omega \times \mathbb{R}}\left[ f^{\delta
}-\left( f^{\delta }\right) ^{2}\right] \ \psi _{2}^{\delta }(\mathbf{x})\
dvd\mathbf{x},
\end{equation*}%
where $f^{\delta }:=\frac{1}{\delta }\int_{0}^{\delta }f(t)dt.$ Since $%
0\leqslant f^{\delta }\leqslant 1,$ in view of the dominated convergence
theorem and Proposition \ref{TRACES}, we derive%
\begin{equation*}
\limsup_{\delta \rightarrow 0}C_{1}^{\delta }\leqslant \int_{\Omega \times
\mathbb{R}}(f^{0}-f{^{0}}^{2})\ dvd\mathbf{x}=0.
\end{equation*}

Let us now consider the term $C_{2}^{\delta }.$ Since $\ \Gamma \in C^{2},$
there exists a small $\delta ,$ such that any point $\mathbf{x}\in S_{\delta
}:=\left\{ \mathbf{x}\in \Omega :\ d(\mathbf{x})<\delta \right\} $ \ has an
unique projection $\mathbf{x}_{0}=\mathbf{x}_{0}(\mathbf{x})$ on the
boundary $\Gamma .$ In the set $S_{\delta },$ we have that $\nabla d(\mathbf{%
x})=-\mathbf{n}(\mathbf{x}_{0})+O(\delta )$ and the Jacobian of the change
of variables $\mathbf{x}\leftrightarrow (\mathbf{x}_{0},s)$ with $s:=d(%
\mathbf{x})$ is equal to $\frac{D(\mathbf{x})}{D(\mathbf{x}_{0},s)}%
=1+O(\delta ),$ since $(\mathbf{x}_{0},s)$ forms the orthogonal coordinate
system at $s=0.$ In view of $\mathbf{v}(t,\cdot )\in \mathbf{V}^{1}(\Omega )$
for a. a. $t\in (0,T),$ we can apply theorems 6.5.3-6.5.4 of \ \cite{Kuf}
and obtain with the help of (\ref{holder}) the following inequality
\begin{equation*}
C_{2}^{\delta }\leqslant \int_{\Gamma _{T}\times \mathbb{R}}|g^{\prime }(v)%
\mathbf{b}_{\mathbf{n}}(t,\mathbf{x}_{0})|^{-}\ [f^{\delta }-\left(
f^{\delta }\right) ^{2}]\ dvdtd\mathbf{x}_{0}+O(\delta ^{\alpha }).
\end{equation*}%
Here $f^{\delta }:=\frac{1}{\delta }\int_{0}^{\delta }f(\cdot ,\cdot ,%
\mathbf{x})ds.$\ Hence Proposition \ref{TRACES} implies%
\begin{equation*}
\limsup_{\delta \rightarrow 0}C_{2}^{\delta }\leqslant \int_{\Gamma
_{T}\times \mathbb{R}}|g^{\prime }(v)\mathbf{b}_{\mathbf{n}}(t,\mathbf{x}%
)|^{-}(f^{b}-\left( f^{b}\right) ^{2})\ dvdtd\mathbf{x}=0.
\end{equation*}%
Finally integrating (\ref{ddF}) over $t_{0}\in \lbrack 2\delta ,T],$
applying Fubini's theorem to the left part of the inequality and taking the
limit on $\delta \rightarrow 0,$ we get \ $\int_{\Omega _{T}\times \mathbb{R}%
}F\ dvdtd{\mathbf{x}}\leqslant 0$. Therefore $F\equiv 0$ a.e. in $\Omega
_{T}\times \mathbb{R}.$ \qquad\ $\blacksquare $

\bigskip

Since $f$ is monotone decreasing on $v$ and $f$ \ takes only the values $0$
and $1,$ \ a. e. in $\Omega _{T}\times \mathbb{R},$ there exists a function $%
z=z(t,\mathbf{x}),$ such that
\begin{equation*}
f(t,\mathbf{x},v)=\mathrm{sgn}^{+}(z(t,\mathbf{x})-v).
\end{equation*}%
Therefore for any $G\in C^{1}([0,1]),G(0)=0$%
\begin{equation*}
\ G\left( u^{\varepsilon }\right) =\int\limits_{0}^{1}G^{\prime
}(v)f^{\varepsilon }(\mathbf{\cdot },\mathbf{\cdot ,}v)\ dv\rightharpoonup
\int\limits_{0}^{1}G^{\prime }(v)f(\mathbf{\cdot },\mathbf{\cdot ,}v)\
dv=G(z)
\end{equation*}%
weakly -- $\ast $ in $L^{\infty }(\Omega _{T}).$ This implies $z=u$ \ and
the strong convergence of $\left\{ u^{\varepsilon }\right\} $\vspace{1pt} to
$u$ \ in $L^{p}(\Omega _{T})$ \ for any $p<\infty .$ Therefore, the function
$\mathbf{v}$ fulfills the integral identity \eqref{V}. And if we take the
sum of the (in)equalities \eqref{EPB+L}, \eqref{EPB-L}, we derive that $u$
satisfies \eqref{DGS}, that ends the proof of Theorem \ref{ETGS}.

\vspace{1pt}

\begin{remark}
\label{rem} Let us note that the measures $m_{+},m_{-}$ and the limit
functions $u,\mathbf{v}$ with $f(t,\mathbf{x},v)=\mathrm{sgn}^{+}(u(t,%
\mathbf{x})-v)$ satisfy all relations \eqref{EPB+L}-\eqref{f} too.
\end{remark}

\vspace{1pt}

\vspace{1pt}

\section{Statement of the quasi-stationary \textbf{Stokes }\textit{B-L}
system}

\label{SSOP1}

For a given viscous parameter $\nu >0,$ we consider the following
initial-boundary value problem, denoted as \textbf{IBVP}$_{\tau =0}$:

\textit{Find a pair }$(u,\mathbf{v})=(u(t,\mathbf{x}),\mathbf{v}(t,\mathbf{x}%
)):\Omega _{T}\rightarrow \mathbb{R}\times \mathbb{R}^{d}$\textit{\ solution
to the quasi-stationary Stokes-Buckley-Leverett system in the domain } $%
\Omega _{T}$
\begin{eqnarray}  \label{EDFV1}
\partial _{t}u+\mathrm{div}\big(\mathbf{v}\;g(u)\big) &=&0,  \label{ETU1} \\
-\nu \Delta \mathbf{v}+h(u)\mathbf{v} &=&-\nabla p,\qquad \mathrm{div}\left(
\mathbf{v}\right) =0,  \label{EPEV1}
\end{eqnarray}%
\textit{satisfying the boundary conditions}%
\begin{equation}
(u,\mathbf{v})=(u_{b},\mathbf{b})\quad \text{on }\Gamma _{T},  \label{BC1}
\end{equation}%
\textit{\ and the initial condition }%
\begin{equation}
u=u_{0}\quad \text{in }\Omega .  \label{IC1}
\end{equation}

\vspace{1pt}

We assume that our data $g,$ $h,\ \ u_{b},$ $u_{0},\ \mathbf{b}$\ satisfy
the following regularity properties%
\begin{eqnarray}
g,\ h &\in &W_{\,\mathrm{loc}}^{1,\infty }(\text{$\mathbb{R}$})\qquad \text{%
with}\quad 0<h_{0}\leqslant h(u),\text{ }  \notag \\
0 &\leqslant &u_{b}\leqslant 1\quad \text{on }\Gamma _{T},  \notag \\
0 &\leqslant &u_{0}\leqslant 1\quad \text{in }\Omega ,  \label{reg33} \\
\mathbf{b} &\in &\mathbf{G}(\Gamma _{T}).  \label{reg3}
\end{eqnarray}

\begin{definition}
\label{DGS12} A pair of functions
\begin{equation*}
(u,\mathbf{v})\in L^{\infty }(\Omega _{T})\times L^{2}(0,T;\mathbf{V}%
^{1}(\Omega ))
\end{equation*}%
is called a weak solution to the \textbf{IBVP}$_{\tau =0}:$ \ \eqref{ETU1}-%
\eqref{IC1}, if the pair $(u,\mathbf{v})$ satisfies the integral inequality%
\begin{align}
\iint_{\Omega _{T}}& \big(|u-v|\;\phi _{t}+\,\mathrm{sgn}(u-v)\big(g(u)-g(v)%
\big)\;\mathbf{v}\cdot \nabla \phi \big)\ d\mathbf{x}\,dt  \notag \\
& +\int_{\Gamma _{T}}M\ |u_{b}-v|\phi \ d\mathbf{x}\,dt+\int_{\Omega
}|u_{0}-v|\phi (0,x)\ d\mathbf{x}\geq 0,  \label{DGV21}
\end{align}%
for any fixed $v\in \mathbb{R},$ where $M:=\emph{K}|\mathbf{b}_{\mathbf{n}}%
\mathbf{|}$ on $\Gamma _{T}$ with $\ \emph{K:=}||g^{\prime }||_{L^{\infty }(%
\mathbb{R})}$\ and for any nonnegative function $\phi \in C_{0}^{\infty
}((-\infty ,T)\times \mathbb{R}^{d})$ and also the following integral
identity
\begin{equation}
\int_{\Omega }\left[ \nu \,\nabla \mathbf{v}:\nabla \boldsymbol{\psi }+h(u)\
\mathbf{v}\cdot \boldsymbol{\psi }\right] \,d\mathbf{x}=0\quad \text{for a.
a. }t\in (0,T)  \label{V21}
\end{equation}%
holds for any $\boldsymbol{\psi }\in \mathbf{C}_{0}^{1}(\Omega ).$ Moreover
the trace of $\mathbf{v}$ is equal to $\mathbf{b}$ on $\Gamma _{T}.$
\end{definition}

\begin{theorem}
\label{ETGS1} If the data $g,$ $h,$ $u_{b},$ $u_{0},$ $\mathbf{b}$ fulfills
the regularity properties \eqref{reg33}-\eqref{reg3}, then the \textbf{IBVP}$%
_{\tau =0}$ \ has a weak solution $(u,\mathbf{v}),$ satisfying%
\begin{eqnarray*}
0 &\leqslant &u\leqslant 1\text{\qquad a . e. in \quad }\Omega _{T}, \\
\mathbf{v,\ }\partial _{t}\mathbf{v} &\in &L^{2}(0,T;\mathbf{V}^{1}(\Omega
)).
\end{eqnarray*}
\end{theorem}

\bigskip

\vspace{1pt}

\subsection{Existence of weak solution. The limit transition on $\protect%
\tau \rightarrow 0$}

\label{SEGS1}

Let us choose some function $\mathbf{v}_{0}\in \mathbf{V}^{0}(\Omega ),$
such that
\begin{equation*}
\mathbf{v}_{0}\cdot \mathbf{n=b}(0)\cdot \mathbf{n}\quad \quad \text{in}%
\quad H^{-1/2}(\Gamma ).
\end{equation*}%
Then, due to Theorem \ref{ETGS}, for each $\tau >0$, there exists a solution
$(u^{\tau },\mathbf{v}^{\tau })$ for the problem \textbf{IBVP}$_{\tau }$: %
\eqref{ETU}-\eqref{IC}, satisfying \eqref{r}. Hereupon the issue is to pass
to the limit on the parameter\ $\tau \rightarrow 0$ \ and, as a consequence,
to derive the solvability of \textbf{IBVP}$_{\tau =0}$.

\begin{proposition}
\label{PPA2} There exists a pair $(u,\mathbf{v})\in L^{\infty }(\Omega
_{T})\times L^{2}(0,T;\mathbf{V}^{1}(\Omega )),$ with $\partial _{t}\mathbf{v%
}\in L^{2}(0,T;\mathbf{V}^{1}(\Omega )),$ and a subsequence of $\left\{
u^{\tau },\mathbf{v}^{\tau }\right\} _{\tau >0},$ such that
\begin{equation}
u^{\tau }\rightharpoonup u\quad \ast \text{-weakly in }L^{\infty }(\Omega
_{T}),  \label{aq1}
\end{equation}%
\begin{equation}
\mathbf{v}^{\tau }\rightarrow \mathbf{v}\qquad \text{strongly in \textbf{$L$}%
$^{2}(\Omega _{T})$}.  \label{aq2}
\end{equation}
\end{proposition}

\textbf{Proof.}\ \ The convergence \eqref{aq1} \ follows from the first
estimate of \eqref{r}. Hence it remains to show \eqref{aq2}.

For each $\tau >0$, let us consider the quasi-stationary Stokes type system
\begin{equation*}
\begin{cases}
-\nu \,\Delta \mathbf{B}^{\tau }+h(u^{\tau })\,\mathbf{B}^{\tau }=-\nabla
\pi ^{\tau },\qquad \mathrm{div}(\mathbf{B}^{\tau })=0\quad \text{in }\Omega
_{T}, \\
\mathbf{B}^{\tau }=\mathbf{b}\quad \text{on }\Gamma _{T}.%
\end{cases}%
\end{equation*}%
The function $\mathbf{z}^{\tau }:=\mathbf{B}^{\tau }-\mathbf{v}_{b},$ \
where $\mathbf{v}_{b}$ is the solution of \eqref{v1}, fulfills the system
\begin{equation*}
\begin{cases}
-\nu \,\Delta \mathbf{z}^{\tau }+h(u^{\tau })\,\mathbf{z}^{\tau }=-\nabla
(\pi ^{\tau }-p_{b})+\mathbf{f}^{\tau },\qquad \mathrm{div}(\mathbf{z}^{\tau
})=0\quad \text{in }\Omega _{T}, \\
\mathbf{z}^{\tau }=\mathbf{0}\quad \text{on }\Gamma _{T},%
\end{cases}%
\end{equation*}%
with $\mathbf{f^{\tau }:=}-h(u^{\tau })\mathbf{v}_{b}.$ Therefore, for $a.a.
\, t \in (0,T)$, $z^\tau(t)$ satisfies the following estimate %
\begin{equation}
\nu \Vert \nabla \mathbf{z}^{\tau }\Vert _{\mathbf{L}^{2}(\Omega
)}^{2}\leqslant \int_{\Omega }|(\mathbf{f}^{\tau }\cdot \mathbf{z}^{\tau
})|\,d\mathbf{x}\leqslant \frac{\nu }{2}\Vert \nabla \mathbf{z}^{\tau }\Vert
_{\mathbf{L}^{2}(\Omega )}^{2}+C\Vert \mathbf{v}_{b}\Vert _{\mathbf{L}%
^{2}(\Omega )}^{2}.  \label{sss}
\end{equation}%
%
%
%
%
%
Now, for $a.a. \, t_0,t_1 \in (0.T)$, we could write %
\begin{equation*}
\ ||\mathbf{v}_{b}\,(t_{1},\cdot )||^{2}-|| \mathbf{v}_{b}(t_{0},\cdot
)\,||^{2}\mathbf{=}\int_{t_{0}}^{t_{1}}\left( \partial _{t}\int_{\Omega }%
\mathbf{v}_{b}^{2}\,d\mathbf{x}\right) \,dt=: J,
\end{equation*}
hence by \eqref{v2} 
\begin{equation*}
|J|\leqslant \Vert \mathbf{v} _{b}\Vert _{\mathbf{L}^{2}(\Omega
_{T})}^{2}+\Vert \partial _{t}\mathbf{v} _{b}\Vert _{\mathbf{L}^{2}(\Omega
_{T})}^{2}\leqslant C
\end{equation*}
and we have that $\mathbf{v}_{b}\in C([0,T];\mathbf{V}^{0}(\Omega _{T}))$.
Consequently, by \eqref{sss}, it follows that
\begin{equation}
\Vert \mathbf{B}^{\tau }\Vert _{L^{\infty }(0,T;\mathbf{V}^{1}(\Omega
))}\leqslant C.  \label{s1}
\end{equation}%
%
%
%
%
Here and below, $C$ are denoted constants, which could change from one to
another statement, being independent of the parameter $\tau $.

Since the function $u^{\tau }$ is the solution of \eqref{ETU} (in the weak
form), the pair $\mathbf{Z}^{\tau }:=\partial _{t}\mathbf{z}^{\tau }$
fulfills the system
\begin{equation}
\begin{cases}
-\nu \,\Delta \mathbf{Z}^{\tau }+h(u^{\tau })\,\mathbf{Z}^{\tau }=-\nabla
Q^{\tau }+\mathbf{R}^{\tau },\qquad \mathrm{div}(\mathbf{Z}^{\tau })=0\quad
\text{in }\Omega _{T}, \\
\mathbf{Z}^{\tau }=0\quad \text{on }\Gamma _{T},%
\end{cases}
\label{s2}
\end{equation}%
with $Q^{\tau }:=\partial _{t}(\pi ^{\tau }-p_{b})$ and $\mathbf{R}^{\tau }$
is given by
\begin{equation*}
\mathbf{R}^{\tau }:= \mathrm{div}\Big(r(u^{\tau }) \big(\mathbf{B}^{\tau }
\otimes \mathbf{v}^{\tau }\big) \Big) - r(u^{\tau}) \Big(\nabla \mathbf{B}%
^{\tau }\Big) \mathbf{v}^{\tau } - h(u^{\tau }) \; \partial _{t}\mathbf{v}%
_{b},
\end{equation*}
%
%
where $r(u):=\int_{0}^{u}h^{\prime }(s)g^{\prime }(s)\ ds.$ From \eqref{s2},
we obtain
\begin{eqnarray*}
\nu \Vert \nabla \mathbf{Z}^{\tau }\Vert _{\mathbf{L}^{2}(\Omega )}^{2}
&\leqslant &\int_{\Omega }(\mathbf{R}^{\tau }\cdot \mathbf{Z}^{\tau })\,d%
\mathbf{x}\leqslant C\Vert \nabla \mathbf{Z}^{\tau }\Vert _{\mathbf{L}%
^{2}(\Omega )} \\
&&\times \left\{ \Vert \mathbf{v}^{\tau }\Vert _{\mathbf{L}^{4}(\Omega )}%
\left[ \Vert \mathbf{B}^{\tau }\Vert _{\mathbf{L}^{4}(\Omega )}+\Vert \nabla
\mathbf{B}^{\tau }\Vert _{\mathbf{L}^{2}(\Omega )}\right] +\Vert \partial
_{t}\mathbf{v}_{b}\Vert _{\mathbf{L}^{2}(\Omega )}\right\} .
\end{eqnarray*}%
The embedding theorem $H^{1}(\Omega )\hookrightarrow L^{4}(\Omega )$\ and %
\eqref{r}, \eqref{v2}, \eqref{s1} imply
\begin{equation}
\Vert \partial _{t}\mathbf{B}^{\tau }\Vert _{L^{2}(0,T;\mathbf{H}^{1}(\Omega
))}\leqslant C.  \label{s5}
\end{equation}

Finally we consider the difference $\mathbf{w}^{\tau }:=\mathbf{v}^{\tau }-%
\mathbf{B^{\tau },}$\ which\ satisfies the system%
\begin{eqnarray*}
\tau \,\partial _{t}\mathbf{w}^{\tau }-\nu \Delta \mathbf{w}^{\tau
}+h(u^{\tau })\mathbf{w}^{\tau } &=&-\nabla (p^{\tau }-\pi \mathbf{^{\tau }}%
)+\tau \mathbf{f}^{\tau },\qquad \mathrm{div}\left( \mathbf{w}^{\tau
}\right) =0, \\
\mathbf{w}^{\tau }\big|_{\Gamma _{T}} &=&0,\qquad \mathbf{w}^{\tau }\big|%
_{t=0}=\mathbf{v}_{0}-\mathbf{B}^{\tau }\big|_{t=0},
\end{eqnarray*}%
with $\mathbf{f^{\tau }:=-}\partial _{t}\mathbf{B}^{\tau }.$ \ \ If we
multiply the first equation of this system by $\mathbf{w}^{\tau }$ and
integrate over $\Omega ,$ we obtain
\begin{align*}
\frac{d}{dt}\left( \frac{\tau }{2}\Vert \mathbf{w}^{\tau }\Vert _{\mathbf{L}%
^{2}(\Omega ))}^{2}\right) +\nu \Vert \nabla \mathbf{w}^{\tau }\Vert _{%
\mathbf{L}^{2}(\Omega )}^{2}& \leqslant C\int_{\Omega }|(\tau \mathbf{f}%
^{\tau }\cdot \mathbf{w}^{\tau })|\,d\mathbf{x} \\
& \leqslant \frac{\nu }{2}\Vert \nabla \mathbf{w}^{\tau }\Vert _{\mathbf{L}%
^{2}(\Omega )}^{2}+C\tau ^{2}\Vert \mathbf{f}^{\tau }\Vert _{\mathbf{L}%
^{2}(\Omega )}^{2}.
\end{align*}%
Integrating the last inequality over the time interval $(0,t)$\ and using (%
\ref{reg3}), (\ref{s5}), we deduce%
\begin{equation}
\Vert \mathbf{v}^{\tau }-\mathbf{B^{\tau }}\Vert _{L^{2}(0,T;\mathbf{V}%
^{1}(\Omega ))}^{2}\leqslant C\tau .  \label{s6}
\end{equation}

\bigskip

Obviously the derived estimates (\ref{s1}), (\ref{s5}), (\ref{s6}) imply the
existence of a function $\mathbf{v\in }L^{2}(0,T;\mathbf{V}^{1}(\Omega )),$
satisfying the strong convergence (\ref{aq2}) for some subsequence of $%
\left\{ \mathbf{v}^{\tau }\right\} _{\tau >0}.\qquad \blacksquare $

\bigskip

Of course, the\ convergence \eqref{aq1}-\eqref{aq2} is not sufficient to
take the limit transition on $\tau \rightarrow 0$ in the system \eqref{ETU1}%
--\eqref{EPEV1}, since we need the strong convergence of a subsequence for $%
\{u^{\tau }\}_{\tau >0}$. To get this strong convergence, we can apply the
Kinetic approach, developed in Section \ref{SEGS} and prove Theorem \ref%
{ETGS1}. In fact, we have to repeat all considerations of the section \ref%
{SD} (see also Remark \ref{rem}), considering the parameter $\tau $, instead
of $\varepsilon $ \ in \eqref{EPB+}--\eqref{m-property2} (without viscous
terms depending on $\varepsilon ).$ \

\bigskip

\section*{Acknowledgements}


The second author were partially supported by FAPERJ through the grant E-26/
111.564/2008 entitled \textsl{"Analysis, Geometry and Applications"}, and by
Pronex-FAPERJ through the grant E-26/ 110.560/2010 entitled \textsl{%
"Nonlinear Partial Differential Equations"}.

\bigskip



\begin{thebibliography}{99}
\bibitem{ant} \textsc{Antontsev S.N., Kazikhov A.V., Monakhov V.N.,} \emph{%
Boundary-value problems in mechanics of non-homogeneous fluids.} Studies in
Math. and its Appl., Vol. 22, North-Holland, 1990.

\bibitem{arb} \textsc{Arbogast T.,} \emph{The existence of weak solutions to
single porosity and simple dual-porosity models of two-phase incompressible
flow,} Nonlinear Anal., \textbf{19} (11) (1992), 1009-1031.

\bibitem{cat} \textsc{Cattabriga L.}, \emph{Su un problema al contorno
relativo al sistema di equazioni di Stokes (Italian),} Rend. Sem. Mat. Univ.
Padova 31, (1961), 308--340.

\bibitem{CF} \textsc{Chen G.-Q., Frid H.,} \emph{Divergence measure fields
and hyperbolic conservation laws.} Arch. Rational Mech. Anal. 147 (1999)
89--118.

\bibitem{chen} \textsc{Chen Z.,} \emph{Degenerate Two-Phase Incompressible
Flow: I. Existence, Uniqueness and Regularity of a Weak Solution,} J. Dif.
Equations, \textbf{171}, Issue 2, (2001), 203-232.


\bibitem{DCFGRO} \textsc{Córdoba D., Gancedo F., Orive R.} \emph{\
Analytical behavior of two-dimensional incompressible flow in porous media.}
\textit{\ J. Math. Physics,} \textbf{48}(6) 065206 (2007);
doi:10.1063/1.2404593 (19 pages).

\bibitem{Dafermos} \textsc{Dafermos C.M., } \emph{Hyperbolic conservation
laws in continuum physics,} 2nd edition. Springer Verlag, 2005.

\bibitem{lellis} \textsc{De Lellis C.,} \emph{Ordinary differential
equations with rough coefficients and the renormalization theorem of
Ambrosio,} Bourbaki Seminar, Preprint, (2007) 1-26.

\bibitem{DiP-L} \textsc{DiPerna R.J., Lions P.L.,} \emph{Ordinary
differential equations, transport theory and Sobolev spaces.} \textit{\
Invent. Math.} \textbf{98}, (1989) 511--547.

\bibitem{Frid} \textsc{Frid H.,} Solution to the Initial Boundary-Value
Problem for the Regularized Buckley-Leverett System, Acta Applicandae
Mathematicae, 38, 239--265 (1995).

\bibitem{galdi} \textsc{Farwig R., Galdi G.P., Sohr H.,} \emph{A New Class
of Weak Solutions of the Navier--Stokes Equations with Nonhomogeneous Data,}
J. Math. Fluid Mech., \textbf{8} (2006) 423--444.

\bibitem{hornung} \textsc{Hornung U.,} \emph{Homogenization and Porous Media,%
} Interdisciplinary Applied Math., Vol. 6, Springer, (1996).

\bibitem{Kuf} \textsc{Kufner A., Jonh O., Fu$\check{c}$ik S.,} \emph{%
Function Spaces.} Noordholf Intern. Publishing, Leyden (1977).

\bibitem{Lad69} \textsc{Ladyzhenskaya O.A., }\emph{The Mathematical Theory
of Viscous Incompressible Flow.} Gordon and Breach, New York-London, 1969.

\bibitem{LSU68} \textsc{Ladyzhenskaya O.A., Solonnikov V.A., Ural'tseva
N.N., }\emph{Linear and quasilinear equations of parabolic type.} American
Mathematical Society, Providence RJ (1968).

\bibitem{len} \textsc{Lenzinger M., Schweizer B.,} \emph{Two-phase flow
equations with outflow boundary conditions in the hydrophobic hydrophilic
case,} Nonlinear Analysis: Theory, Methods \& Applications, \textbf{73},
Issue 4 (2010), 840-853.

\bibitem{LPT1} \textsc{Lions P.-L., \ Perthame B., \ Tadmor E.,} \emph{%
Kinetic formulation for isentropic gas dynamcs and }$p$\emph{-systems.}
Comm. Math. Phys. 163 (1994), 415--431.

\bibitem{LPT2} \textsc{Lions P.-L., \ Perthame B., \ Tadmor E.,} \emph{A
kinetic formulation of multidimensional scalar conservation laws and related
equations.} J. AMS 7 (1994), 169--191.

\bibitem{SLPIP} \textsc{Luckhaus S., Plotnikov P.I.,} \emph{\ Entropy
solutions to the Buckley-Leverett equations.} Siberian Mathematical Journal
41, N. 2 (2000), 169--191.

\bibitem{malek} \textsc{Malek J., Necas J., Rokyta M., Ruzicka M.,} \emph{%
Weak and measure-valued solutions to evolutionary PDEs.} Chapman\&Hall,
London (1996).

\bibitem{WN1} \textsc{Neves W.,} \emph{Scalar multidimensional conservation
laws IBVP in noncylindrical Lipschitz domains,} Journal of Diff. Equations
192 (2003) 360--395.

\bibitem{O} \textsc{Otto F.,} \emph{Initial-boundary value problem for a
scalar conservation law,} C.R. Acad. Sci. Paris 322 (1996) 729--734.

\bibitem{EP} \textsc{Panov E. Yu.,} \emph{Existence and strong
pre-compactness properties for entropy solutions of a first-order
quasilinear equation with discontinuous flux,} Arch. Ration. Mech. Anal.,
\textbf{195}, no 2, 643--673 (2010).

\bibitem{MPVS} \textsc{Perepetlitsa I., Shelukhin V.,} \emph{On Global
Solutions of a Boundary-Value Problem for the one-~imensionaBl
uckley-Leverett Equations,} Applicable Analysis, \textbf{73}, no 3--4,
325--343 (1999).

\bibitem{PD} \textsc{Perthame B., Dalibard A.-L.,} \emph{Existence of
solutions of the hyperbolic Keller-Segel model, }Trans. Amer. Math. Soc.,
\textbf{361}, 2319-2335 (2009).

\bibitem{perthame2} \textsc{Perthame B., }\emph{Kinetic formulation of
conservation laws, }Oxford University Press, 2002.

\bibitem{SAS} \textsc{Sazhenkov S. A.,} \emph{\ Entropy solutions to the
Verigin ultraparabolic problem,} Siberian Mathematical Journal 49, No. 2,
362--374 (2008).

\bibitem{AES} \textsc{Scheidegger A.E.}, \emph{Hydrodynamics in Porous Media}%
, Handbuch der Physik Vol. VIII/2, Flûgge, Springer, (1963).

\bibitem{AES1} \textsc{Scheidegger A.E.,} \emph{The Physics of Flow Through
Porous Media}, 3rd ed, University of Toronto Press, Toronto (1974).

\bibitem{sheu} \textsc{Sheu L.-J.,} \emph{An autonomous system for chaotic
convection in a porous medium using a thermal non-equilibrium model,} Chaos,
Solitons and Fractals, \textbf{30} (2006) 672--689.

\bibitem{straughan} \textsc{Straughan B.,} \emph{Stability and Wave Motion
in porous media,} Applied Math. Sciences Vol. 165, Springer, (2008).

\bibitem{Temam} \textsc{Temam R.,} \emph{Navier-Stokes equations,} Theory
and numerical analysis. AMS Chelsea publishing, Providence, Rhode Island
(2001).

\bibitem{wang} \textsc{Wang B., Lin S.,} \emph{Existence of global
attractors for the three-dimensional Brinkman Forchheimer equation, }Math.
Meth. Appl. Sci., \textbf{31} (2008), 1479--1495.
\end{thebibliography}
\end{document}